\documentclass{amsart}

\usepackage{graphicx}
\usepackage{color}
\usepackage{amsmath}
\usepackage{amsfonts}
\usepackage{amsopn}
\usepackage{amssymb}

\newcommand{\sketch}{\noindent \emph{Sketch of proof: }}
\newcommand{\R}{\ensuremath{\mathbb{R}}}

\newcommand{\Z}{\ensuremath{\mathbb{Z}}}
\newcommand{\C}{\ensuremath{\mathbb{C}}}
\newcommand{\Quat}{\ensuremath{\mathbb{H}}}
\newcommand{\Octon}{\ensuremath{\mathbb{O}}}
\newcommand{\K}{\ensuremath{\mathbb{K}}}
\newcommand{\Heis}{\ensuremath{H}}
\newcommand{\bd}{\ensuremath{\partial}}

\newcommand{\Clip}{C^{\text{lip}}}

\DeclareMathOperator{\interior}{int}
\DeclareMathOperator{\area}{area}
\DeclareMathOperator{\supp}{supp}
\DeclareMathOperator{\Lip}{Lip}
\DeclareMathOperator{\id}{id}

\DeclareMathOperator{\vol}{vol}
\DeclareMathOperator{\FV}{FV}
\DeclareMathOperator{\mass}{mass}
\renewcommand{\setminus}{{-}}
\def\cprime{$'$}

\newtheorem{thm}{Theorem}
\newtheorem{dfn}{Definition}
\newtheorem{lemma}{Lemma}[section]
\newtheorem{prop}[lemma]{Proposition}
\newtheorem{cor}[lemma]{Corollary}

\title{Filling inequalities for nilpotent groups}
\author{Robert Young}
\address{New York University\\
251 Mercer St., New York, NY  10012}
\date{\today}
\email{rjyoung1729@gmail.com}

\begin{document}
\bibliographystyle{hamsplain}

\begin{abstract}  
  We bound the higher-order Dehn functions and other filling
  invariants of certain Carnot groups using approximation techniques.
  These groups include the higher-dimensional Heisenberg groups, jet
  groups, and central products of two-step nilpotent groups.  Some
  consequences of this work are a construction of groups with
  arbitrarily large nilpotency class that have euclidean
  $n$-dimensional filling volume functions, and a proof of part of a
  conjecture of Gromov on the higher-order filling functions of the
  higher-dimensional Heisenberg groups.
\end{abstract}

\maketitle 

\section{Introduction}
A filling invariant is an invariant of a space which measures the
difficulty of finding a disc or chain with a specified boundary.
These invariants are studied in geometric group theory because their
asymptotic growth rates describe the large-scale geometry of a group
or a space.

The primary example of a filling invariant is the Dehn function of a
space, which describes the area necessary to fill a closed curve with
a disc.  This is an important invariant of geometric group theory
because if $G$ is a group that acts geometrically (cocompactly,
properly discontinously, and by isometries) on a space, then filling a
closed curve in the space with a disc corresponds to reducing a word
in the group to the identity using relators.  Thus, the Dehn function
gives a way to connect the combinatorial group theory of a group to
the geometry of spaces that it acts on.

The Dehn function can be generalized to the $k$th-order Dehn function
$\delta_X^k(V)$ of a space $X$, which bounds the volume necessary to
fill a $k$-sphere in $X$ of volume $\le V$ by a $k+1$-disc.  Since spheres
and cycles can have much more complicated metrics than closed curves,
much less is known about these higher-order filling invariants than
about Dehn functions.  In this paper, we will give methods for
bounding the higher-order Dehn functions of certain nilpotent groups
based on techniques from geometric measure theory.

These techniques generalize work of Gromov on the Dehn functions of
nilpotent groups.  Gromov proved that Lipschitz maps of spheres into
nilpotent groups which have sufficiently many horizontal maps (i.e.,
groups which satisfy an appropriate flexibility condition) can be
extended to Lipschitz maps of balls \cite{GroCC}.  In particular, if a
group satisfies the conditions, a closed curve in the group of length
$\ell$ can be filled by a Lipschitz map of a disc with Lipschitz
constant $O(\ell)$, providing a Dehn function estimate \cite{GroAII};
this works because a closed curve can be replaced by a Lipschitz map
of a circle with Lipschitz constant $\ell$.  These techniques do not,
however, give good higher-order Dehn function estimates, since a map
of a sphere with bounded volume may have an arbitrarily large
Lipschitz constant.

We will solve this problem by using approximation techniques and
combining information from multiple scales.  We use the
Federer-Fleming Deformation Theorem to approximate singular cycles in
a group by simplicial cycles; by rescaling the cycles, we produce
approximations which use larger or smaller simplices.  By linking
these approximations together, we get a filling of the original
cycle.  Under appropriate conditions, the approximations
and the filling all have small volumes.

This method constructs a filling of a sphere out of scalings of a
finite set of pieces; indeed, from scalings of the simplices of a
triangulation satisfying certain properties.  Any way of constructing
such a triangulation gives a bound on the Dehn functions or
higher-order Dehn functions.  We will construct these triangulations
in two different ways: we will use holonomic maps in jet spaces to
construct fillings in jet groups, and we will show that certain
presentations of nilpotent groups can be used to prove quadratic Dehn
functions.

We prove bounds for two families of groups.  The first family consists
of the jet groups.  These groups have been
used as examples of non-rigid Carnot groups by Warhurst \cite{War} and
as a family of quadratically presented Lie algebras by Chen
\cite{Chen}.  We show that these groups have many horizontal
submanifolds and use these submanifolds to construct triangulations;
this leads to euclidean bounds on their filling functions.
Specifically, we show:
\begin{thm}\label{jdehn}
For any $m>0,k>0$, the group $J^m(\R^k)$, when endowed with a left-invariant riemannian metric, has
$$\FV^i(V)\sim \delta^{i-1}\sim V^{\frac{i}{i-1}}$$
for $i\le k$, and
$$\FV^{k+1}(V)\sim \delta^{k}\sim V^{\frac{k+m+1}{k}}.$$
\end{thm}
Here, $J^m(\R^k)$ is a $(k+1)$-step nilpotent group based on the
$m$-jet bundle of $\R^k$.  For $i\le k$, these are the same filling
inequalities as for euclidean space, so we call these bounds {\em
  euclidean}.  This gives examples of groups of any nilpotency class
with euclidean filling functions.  

This theorem answers open questions about the filling functions of the
Heisenberg groups and the existence of nilpotent groups of large
nilpotency class with quadratic Dehn functions.  Gromov conjectured
that the filling functions of the $(2k+1)$-dimensional Heisenberg
group $\Heis^{2k+1}$ should be euclidean in dimensions below $k$,
super-euclidean in dimension $k$, and sub-euclidean in dimensions
above $k$ \cite[Conj.\ 5.D.5(f)]{GroAII}.  Since
$J^1(\R^k)=\Heis^{2k+1}$, where $\Heis^{2k+1}$ is the
$2k+1$-dimensional Heisenberg group, Theorem~\ref{jdehn} implies two
parts of the conjecture:
\begin{cor}
For $\Heis^{2k+1}$ endowed with a left-invariant riemannian metric,
$$\FV^i(V)\sim \delta^{i-1}\sim V^{\frac{i}{i-1}}$$
for $i\le k$, and
$$\FV^{k+1}(V)\sim \delta^{k}\sim V^{\frac{k+2}{k}}.$$
\end{cor}

Furthermore, the jet groups give the first examples of nilpotent
groups of arbitrary nilpotency class which have quadratic Dehn
function.  This is surprising because nilpotent groups with large
nilpotency class have highly distorted subgroups; in a class $c$
torsion-free nilpotent group, there are elements $x$ such that $x^t$
is distance $\sim t^{1/c}$ away from the identity in the word metric.
As a consequence, the Dehn functions of nilpotent groups with large
nilpotency class are often large; for instance, the free $c$-step
nilpotent groups have Dehn functions $\delta(n)\sim n^{c+1}$.

The second class of groups we consider are nilpotent groups
constructed as quotients by commutators.  An important special case is
the case of central products of nilpotent groups.  These groups
generalize the construction used to construct the higher-dimensional
Heisenberg groups.  Allcock \cite{Allcock} and Ol'shanskii-Sapir
\cite{SapOls} studied the higher-dimensional Heisenberg groups with
methods that generalize to other central products (see
e.g.\ \cite{Magnani} for one generalization).  We will show that
bounds on the Dehn functions of many central products follow from these
quotient theorems.

In Section~\ref{prelims}, we describe some of the objects that we will
use in this paper, such as scaling automorphisms of Carnot groups,
Lipschitz chains, and Lipschitz triangulations.  In
Section~\ref{sec:fillApprox}, we describe how one can use
triangulations of groups and approximations by simplicial chains to
construct fillings of cycles in a nilpotent group.  In order to get
good bounds on filling inequalities, these triangulations must satisfy
certain metric properties, and we construct such triangulations in
Section~\ref{sec:horizTris}.  Finally, in Sections~\ref{sec:moreUpper}
and \ref{sec:lower}, we finish the proof of Theorem~\ref{jdehn} by
showing upper bounds in the middle dimension and by using
cohomological techniques to prove lower bounds on filling inequalities
in jet groups.

\noindent \emph{Note:} Much of the work in this paper was done as part
of my doctoral thesis at the University of Chicago, and many of the
results in this paper first appeared in preprint form as
\cite{YoungSca} and \cite{YoungFING}.  The exposition and proofs,
however, have been substantially rewritten and improved.  The author
would like to thank Shmuel Weinberger, Benson Farb, and Stefan Wenger
for their help and advice, and to thank an anonymous referee for their
advice on clarifying and simplifying the exposition.

\section{Preliminaries}\label{prelims}
One inspiration for the study of Dehn functions is the classical
isoperimetric problem, which asks for the largest area in the plane
which can be enclosed by a loop of a given length.  This can be
generalized to other spaces and other boundaries; instead of asking
about the area enclosed by a loop in the plane, we can ask about the
filling volume of a $k$-sphere or $k$-cycle in some other space (i.e.,
the infimal volume of a $(k+1)$-ball or chain with that sphere or
cycle as a boundary).  This leads to several possible definitions of
higher-order filling invariants, each using a different class of
fillings or of boundaries.  In this paper, we will mainly use a
definition based on Lipschitz chains (i.e., formal sums of Lipschitz
maps of simplices), as in \cite{GroFRM} and \cite{WengerSimp}.  Other
definitions of higher-order filling invariants can be found in
\cite{AlWaPr, ECHLPT, GroftA, GroftB}.

Recall that an \emph{integral Lipschitz singular $d$-chain} in a space $X$ is a finite
linear combination, with integer coefficients, of Lipschitz maps from
the euclidean $d$-dimensional simplex $\Delta^d$ to $X$.  We will
often call this simply a \emph{Lipschitz $d$-chain}.

We will generally take $X$ to be a riemannian manifold or simplicial
complex, so by Rademacher's Theorem, a Lipschitz map from a simplex to
$X$ is differentiable almost everywhere, and one can define the volume
of the map as the integral of the magnitude of its jacobian.  If $a$
is a Lipschitz $d$-chain, for $d>0$ and $a=\sum_i x_i \alpha_i$ for
some maps $\alpha_i:\Delta^d\to X$ and some coefficients $x_i\in \Z,
x_i\ne 0$, we define the {\em mass} of $a$ to be
$$\mass{a}:=\sum_i x_i \vol_d \alpha_i.$$
The Lipschitz chains form a chain complex, which we denote
$\Clip_*(X)$.  If $g:X\to Y$ is a Lipschitz map, we can
define the pushforward map, $g_\sharp:\Clip_*(X)\to \Clip_*(Y)$, to be
the linear map which sends the simplex $\alpha:\Delta^d\to X$ to the
simplex $g\circ \alpha$; this is a map of chain complexes.

If $X$ is a $d$-connected riemannian manifold or simplicial complex
and $a$ is an integral Lipschitz $d$-cycle in $X$, we define the
filling volume of $a$ to be:
$$\FV^{d+1}_X(a):=\inf \{\mass{b}\mid \partial b=a\},$$ 
where the infimum is taken over the set of $b\in
\Clip_{d+1}(X)$ such that $\partial b=a$.  We get the
$(d+1)$-dimensional filling volume function by taking a supremum over
all cycles of a given volume:
$$\FV^{d+1}_X(V):=\sup\{\FV^{d+1}_X(a)\mid \mass a \le V\},$$
where the supremum is taken over integral Lipschitz $d$-cycles of mass
at most $V$.

The $(d+1)$-dimensional filling volume function is related to the
$d$-th order Dehn function, $\delta^d$, which measures the volume
necessary to extend a map $S^d\to X$ to a map $D^{d+1}\to X$.  If $X$
is a $d$-connected manifold or simplicial complex and $f:S^d\to X$ is
a Lipschitz map, we define
$$\delta^d_X(f):=\inf \{ \vol_d{g}\mid g:D^{d+1}\to X, g|_{S^d}=f\}$$
and
$$\delta^d_X(V):=\sup \{\delta^d_X(f)\mid f:S^{d}\to X,\vol_d f\le V
\}.$$
When $d=1$, this is simply called the Dehn function, and it is often
written as simply $\delta_X(V)$ or even $\delta(V)$ when the space is clear.

The exact relationship between $\delta^d$ and $\FV^{d+1}$ depends on
$d$.  When $d\ge 3$, then $\delta^d_X\sim \FV^{d+1}_X$ for all
$d$-connected manifolds or simplicial complexes $X$.  When $d=2$, then
$\delta^d_X\lesssim \FV^{d+1}_X$.  (See \cite[App. 2.(A')]{GroFRM} or
\cite{WhiteMinimal} for the upper bound on $\delta^d$ and
\cite[Rem.~2.6.(4)]{BBFS} for the lower bound; see also \cite{GroftA,
  GroftB}.)  When $d\le 2$, the two functions may differ
\cite{YoungHom}, but the bounds found in this paper will hold for
both.  

A key tool to work with these chains is the Federer-Fleming
Deformation Theorem, which states that Lipschitz chains can be
approximated by simplicial chains.  We will state this theorem in
terms of Lipschitz triangulations.  A \emph{triangulation} of a space
$X$ consists of a simplicial complex $\tau$ and a homeomorphism
$f:\tau\to X$.  We will sometimes refer to the tuple $(\tau,f)$
just as $\tau$, leaving the homeomorphism implicit.  We can put a
metric on $\tau$ so that each simplex is isometric to the standard
euclidean simplex; if $f$ is a bilipschitz map (i.e., $f$ and
$f^{-1}$ are both Lipschitz), then we call $(\tau,f)$ a
\emph{Lipschitz triangulation}.  In this paper, all triangulations
will be Lipschitz triangulations.

If $(\tau,f:\tau\to X)$ is a triangulation of $X$, then a
simplicial $k$-chain of $\tau$ is a formal sum of $k$-dimensional
faces of $\tau$.  We can use the push-forward map $f_\sharp$ to
identify such chains with Lipschitz chains of $X$.  By abuse of
notation, we call these chains simplicial chains in $X$, and we denote
the complex of simplicial chains by $C_*(\tau)$.  Federer and Fleming
showed that Lipschitz chains can be approximated by simplicial chains
\cite{FedFlem}.  The version we state is a slight simplification; it
applies only to Lipschitz chains with simplicial boundaries.  In this
theorem, $P_\tau(\alpha)$ will be a simplicial approximation of
$\alpha$ and $Q_\tau(\alpha)$ will be a Lipschitz chain which
interpolates between $\alpha$ and $P_\tau(\alpha)$.
\begin{thm}[Deformation
  Theorem \cite{FedFlem,ECHLPT}]\label{thm:fedflem}
  Let $(\tau,f:\tau\to X)$ be a triangulation of $X$.
  
  There is a constant $c=c(\tau)$ such that if $\alpha\in
  \Clip_k(X)$ is a chain such that $\partial \alpha\in
  C_{k-1}(\tau)$, then there are $P_\tau(\alpha)\in
  C_k(\tau)$ and $Q_\tau(\alpha)\in \Clip_{k+1}(X)$ such that
  \begin{enumerate}
  \item $\mass P_\tau(\alpha)\le c \cdot\mass(\alpha)$,
  \item $\mass Q_\tau(\alpha)\le c \cdot\mass(\alpha)$, and
  \item $\partial Q_\tau(\alpha) = \alpha - P_\tau(\alpha)$.
  \end{enumerate}
\end{thm}
Federer and Fleming originally proved their theorem in the case of
Lipschitz currents in $\R^n$, so this statement is somewhat different
from their original; it is closest to the version proved in \cite{ECHLPT}.

We have defined $\delta^d$ and $\FV^{d+1}$ as invariants of spaces,
but they can in fact be defined as invariants of groups.  If $G$ is a
group which acts geometrically (that is,
cocompactly, properly discontinuously, and by isometries) on a
$d$-connected manifold or simplicial complex $X$, then the asymptotic
growth rates of $\delta^d_X(V)$ and $\FV^{d+1}_X(V)$ are invariants of
$G$.  To make this rigorous, we can define a partial ordering
on functions $\R^+\to \R^+$ so that $f\precsim g$ if and only if there
is a $c$ such that
$$f(x)\le c g(c x + c) + cx + c.$$
We say $f\sim g$ if and only if $f\precsim g$ and $g \precsim f$.
Then if
$X_1$ and $X_2$ are two $d$-connected manifolds or simplicial complexes on which $G$ acts geometrically, then
$$\delta^d_{X_1}(V)\sim \delta^d_{X_2}(V)$$
and
$$\FV^d_{X_1}(V)\sim \FV^d_{X_2}(V)$$
and we thus define $\delta^d_G(V)$ and $\FV^{d+1}_G(V)$ to be the
$\sim$-equivalence classes of $\delta^d_X(V)$ and $\FV^{d+1}_X(V)$.
This is proved for a simplicial version of $\delta^d$ in
\cite{AlWaPr}, but the proof there also applies to a simplicial
version of $\FV^d$; one can show that the Lipschitz versions used here
are equivalent to simplicial versions using the Deformation Theorem.

In this paper, we focus on finding filling inequalities for Carnot
groups, which are nilpotent Lie groups provided with a family of
scaling automorphisms.  Recall that if $G$ is a simply-connected
nilpotent Lie group and $\mathfrak{g}$ is its Lie algebra, then
the lower central series
$$\mathfrak{g}=\mathfrak{g}_0\supset \dots \supset
\mathfrak{g}_{k-1}=\{0\},$$
$$\mathfrak{g}_{i+1}=[\mathfrak{g}_i,\mathfrak{g}]$$
terminates.  If $\mathfrak{g}_k=\{0\}$ and $\mathfrak{g}_{k-1}\ne \{0\}$,
we say that $\mathfrak{g}$ has nilpotency class $k$.  If there is a
decomposition
$$\mathfrak{g}=V_1\oplus \dots \oplus V_k$$
such that 
$$\mathfrak{g}_i=V_{i+1}\oplus \dots \oplus V_k$$
and $[V_i,V_j]\subset V_{i+j}$ for all $i,j\le k$, we call it a
grading of $\mathfrak{g}$.  If $\mathfrak{g}$ has a grading, we can extend the
$V_i$ to left-invariant plane fields on $G$ and give $G$ a
left-invariant metric such that the $V_i$'s are orthogonal.  With this
metric, $G$ is called a Carnot group.

If $G$ is a Carnot group, there is a family of automorphisms $s_t:G\to
G$ which act on the Lie algebra by $s_t(v)=t^i v$ for all $v\in V_i$.
These automorphisms distort vectors in $\mathfrak{g}$ by differing
amounts.  Vectors in $V_1$ are distorted the least, and we call these
vectors {\em horizontal}.  If $M$ is a manifold and $f:M\to
G$ is a Lipschitz map, it is differentiable almost everywhere by
Rademacher's Theorem.  If all of its tangent vectors lie in the
plane field $V_1$, we say that $f$ is {\em horizontal}; likewise, if
$\alpha\in \Clip(G)$, we say that $\alpha$ is horizontal if it is a
sum of horizontal maps.  If
$\tau$ is a simplicial complex and $f:\tau\to G$ is a Lipschitz map
which is horizontal on every simplex of $\tau$ of dimension at most
$k$, we say that $f$ is \emph{$k$-horizontal}.  We then have the
following:
\begin{lemma}
  If $M$ is a compact $k$-dimensional manifold and $f:\Delta^k\to G$ is
  horizontal, then $\vol(s_t\circ f)=t^k\vol f$ for all $t\ge 0$.
\end{lemma}

\section{Filling cycles through approximations}\label{sec:fillApprox}

\subsection{Sketch of argument}
\begin{figure}
  \begin{center}
    \def\svgwidth{\textwidth}
    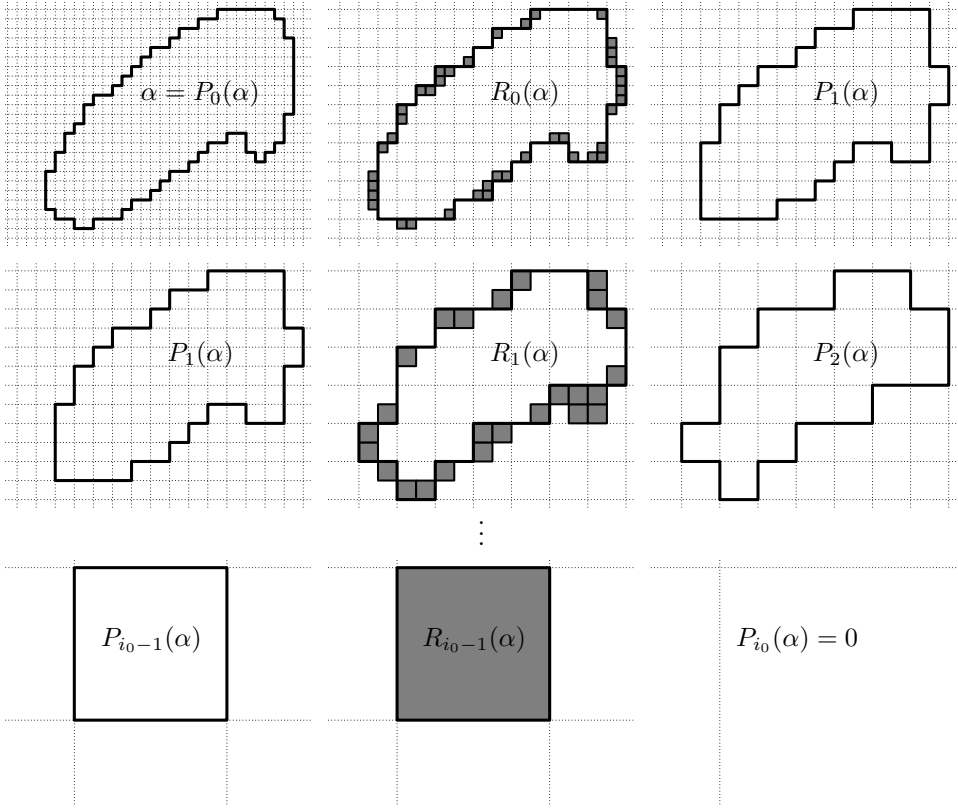
  \end{center}
  \caption{Filling a curve by approximations\label{fig:approxFill}}
\end{figure}
The basic idea behind our techniques is illustrated in Fig.\
\ref{fig:approxFill}, which illustrates a method of bounding the Dehn function
of $\R^2$.  In the figure, $\alpha$ is a curve of length $\ell$ in
$\R^2$.  We approximate $\alpha$ in successively larger grids:
$P_0(\alpha)$ approximates $\alpha$ in a $1\times 1$ grid, and
$P_i(\alpha)$ approximates $\alpha$ in a $2^i\times 2^i$ grid.  Each
approximation has length comparable to the length of $\alpha$ when $i$
is small (i.e., $i\le \log_2\ell$).  When $i> \log_2\ell$, then $\alpha$ is
smaller than any square in the grid and the
curve can be approximated by the zero cycle.  We can
connect approximations using annuli made of squares, and since
$P_i(\alpha)$ and $P_{i+1}(\alpha)$ are close together, it takes
relatively few such squares; just as $P_i(\alpha)$ is made up of $\sim
\ell 2^{-i}$ segments of length $2^i$, $R_i(\alpha)$ is a sum of $\sim
\ell 2^{-i}$ squares with side $2^i$, and thus has area $\sim \ell
2^i$.  If $i_0$ is such that $2^{i_0}\gg \ell$, then 
$P_{i_0}(\alpha)=0$, so we get a filling of $\alpha$ by taking the sum
of the $R_i(\alpha)$'s; if we let
$$\beta=\sum_{i=0}^{i_0} R_i(\alpha),$$
then $\partial \beta=\alpha$ and the area of $\beta$ is $\sim \ell^2$.

A similar argument can be used to fill higher-dimensional cycles in
higher-dimensional euclidean spaces; if $\alpha$ is a $k$-cycle of
mass $V$ in $\R^n$, it can be approximated by a sum $P_i(\alpha)$ of
$\sim V/2^{ik}$ $k$-cubes of side-length $2^i$ by using a $2^i\times
\dots \times 2^i$ grid in $\R^n$.  Furthermore, $P_i(\alpha)$ and
$P_{i+1}(\alpha)$ can be connected by a chain $R_i(\alpha)$ consisting
of $\sim V/2^{ik}$ $(k+1)$-cubes of side-length $2^i$.  As before, if
$i_0$ is such that $2^{i_0k}\gg V$, then $\alpha$ is smaller than any
individual cube, so $P_{i_0}(\alpha)=0$ and if we let
$$\beta=\sum_{i=0}^{i_0} R_i(\alpha),$$
then we get a filling of $\alpha$ with volume $\sim
\ell^{\frac{k+1}{k}}$.

The goal of this paper is to extend this argument to certain nilpotent
groups.  The key step in doing so is to construct the $P_i(\alpha)$
and $R_i(\alpha)$; i.e., to construct a sequence of coarser and
coarser approximations and then to connect those approximations by
chains.  In euclidean space, we use the Federer-Fleming Deformation
Theorem to construct $P_i$; this theorem allows us to approximate
Lipschitz chains by simplicial chains of similar mass, so if we apply
it to a $2^i\times \dots \times 2^i$ grid, we get an approximation by
cubes of side-length $2^i$.  The connecting chains $R_i$ are also
produced by simplicial approximation, but the construction is more
involved.

The bounds we get on the filling functions of nilpotent groups thus
depend on how efficiently we can produce simplicial approximations.
The mass of a simplicial approximation is governed by the constant
factor $c$ in the Federer-Fleming Deformation Theorem, which depends
on the triangulation used.  In $\R^n$, scaling maps allow us to
construct a family of triangulations with differently-sized simplices
and the same $c$, but this is not always possible in a nilpotent
group.  In the rest of this section, we will show that if $G$ is a
Carnot group and sufficiently many horizontal maps into $G$ exist,
then we can produce efficient simplicial approximations at all scales
and use these to find strong bounds on the filling functions of $G$.

Note that we will state and prove our theorems primarily in terms of
$\FV^{k+1}$ rather than $\delta^k$.  At the end of the section, we
will outline how to modify these arguments to provide bounds on
Dehn functions and their higher-order counterparts.

\subsection{Simplicial approximations in Carnot groups}

We will construct $P_i(\alpha)$ and $R_i(\alpha)$ using the
Federer-Fleming Deformation Theorem (Thm.\ \ref{thm:fedflem}).  This
is a tool of geometric measure theory which approximates Lipschitz
chains and cycles in $G$ by simplicial chains and cycles in a
triangulation of $G$.  When $G$ is Carnot, we can construct 
triangulations of $G$ by scaling a single triangulation $\tau$; this gives
triangulations with simplices of different scales, and approximations
in these triangulations give the $P_i$ and $R_i$.  

One difficulty is that the scaling automorphism may distort these
triangulations.  The scaling automorphism $s_t:\R^n\to \R^n$ stretches
each direction by a factor of $t$.  In a Carnot group, however, the
scaling automorphism $s_t:G\to G$ may scale vectors by up to $t^c$,
where $c$ is the nilpotency class of $G$.  Since the $P_i$ and $R_i$
are made up of scaled simplices, a bad choice of $\tau$ may lead to
very large approximations.  We avoid this by using $k$-horizontal
triangulations.  Specifically, we will show that if certain
$k$-horizontal maps and triangulations exist, then $P_i(\alpha)$ can
be constructed so that $\mass P_i(\alpha)\preceq \mass \alpha$ for all
$i$, and that this leads to bounds on the filling volume function.
(In the next section, we will construct some groups which
have such triangulations.)

In the rest of this paper, $G$ will represent a Carnot group and
$\Gamma$ a lattice in $G$.  The maps $s_t:G\to G$, $t\ge0$, will
represent the family of scaling automorphisms of $G$.

Let $(\tau,f:\tau\to G)$ be a triangulation of $G$ and let $P_\tau$
and $Q_\tau$ be as in Theorem~\ref{thm:fedflem}.  If $f$ is
$k$-horizontal and $\alpha$ is a Lipschitz $k$-cycle, then
$P_\tau(\alpha)$ is a horizontal approximation of $\alpha$.  In order
to avoid possible issues with constructing horizontal triangulations
of $G$, we will construct a slightly different approximation.  Let
$(\tau,f:\tau\to G)$ be a triangulation of $G$, with no
conditions on the horizontality of $f$.  Let $k>0$ be an integer and
let $\phi:\tau \to G$ be a $k$-horizontal Lipschitz map which is a
bounded distance from $f$, i.e., there is a $c$ such that
$d(f(x),\phi(x))<c$ for all $x\in G$.  If $\alpha$ is a Lipschitz
$k$-cycle in $G$, then $f^{-1}_\sharp (\alpha)$ is a Lipschitz
$k$-cycle in $\tau$, and $P_\tau(f^{-1}_\sharp (\alpha))$ is a
simplicial cycle in $\tau$.  By abuse of notation we define
$P_{\phi(\tau)}(\alpha)$ by
$$P_{\phi(\tau)}(\alpha):=\phi_\sharp[P_\tau(f^{-1}_\sharp(\alpha))].$$
This is a sum of images of $k$-simplices of $\tau$, so it is a
horizontal cycle.

We will show that the cycle $P_{\phi(\tau)}(\alpha)$ is close to $\alpha$ in
the sense that there is a chain $Q_{\phi(\tau)}(\alpha)$ with mass
comparable to $\mass \alpha$ which interpolates between $\alpha$ and
$P_{\phi(\tau)}(\alpha)$.  
\begin{lemma}\label{lem:Qprops}
  There is a $c_Q$ depending only on $\phi$, $\tau$, and $k$ such that
  for all $k$-cycles $\alpha$, there is a $(k+1)$-chain
  $Q_{\phi(\tau)}(\alpha)$ such that
  $$\partial Q_{\phi(\tau)}(\alpha)=P_{\phi(\tau)}(\alpha)-\alpha$$
  and
  $\mass Q_{\phi(\tau)}(\alpha)\le c_Q \mass \alpha$.
\end{lemma}
\begin{proof}
  Because $\phi\circ f^{-1}$ moves each point of $G$ by a bounded
  distance, there is a Lipschitz homotopy $h:G\times [0,1]\to G$
  between $\id_G$ and $\phi \circ f^{-1}$.  Let
  $$Q_{\phi(\tau)}(\alpha):=h_\sharp(\alpha\times[0,1])+\phi_\sharp[Q_\tau(f^{-1}_\sharp(\alpha))].$$

  For the first part of the lemma, note that
  \begin{align*}
    \partial Q_{\phi(\tau)}(\alpha) &= h_\sharp(\partial(\alpha\times[0,1]))+\phi_\sharp[\partial
    Q_\tau(f^{-1}_\sharp(\alpha))]. \\
    &=(\phi_\sharp(f^{-1}_\sharp(\alpha))-\alpha)+\phi_\sharp[P_\tau(f^{-1}_\sharp(\alpha))-f^{-1}_\sharp(\alpha)] \\
    &=\phi_\sharp(f^{-1}_\sharp(\alpha))-\alpha+\phi_\sharp[P_\tau(f^{-1}_\sharp(\alpha))]-\phi_\sharp(f^{-1}_\sharp(\alpha))=P_{\phi(\tau)}(\alpha)-\alpha.
  \end{align*}
  
  By Theorem~\ref{thm:fedflem}, there is a $c$ such that 
  $$\mass  Q_\tau(f^{-1}_\sharp(\alpha))\le c \mass f^{-1}_\sharp(\alpha).$$
  If we let
  $$c_Q:=\Lip(h)^{k+1}+\Lip(f^{-1})^k\Lip(\phi)^{k+1}c,$$
  it is then straightforward to check the second part of the lemma.
\end{proof}

We can compose $P_{\phi(\tau)}$ with scaling automorphisms to produce
a sequence of approximations.  To avoid cumbersome
subscripts, we will abuse notation by writing $s_{t}(\alpha)$ instead
of $(s_{t})_\sharp(\alpha)$ when the intention is clear.  We define
$$P_i(\alpha)=s_{2^i}( P_{\phi(\tau)}(s_{2^{-i}}(\alpha)));$$
this is a horizontal approximation of $\alpha$ with simplices of diameter
$\sim 2^i$.
\begin{lemma}\label{lem:Pprops}
  There is a $c_P$ depending only on $\phi$ and $\tau$ such that for all
  $i\ge 0$ and for all $k$-cycles $\alpha$,
  we have $\mass P_i(\alpha)\le c_P \mass \alpha$.
\end{lemma}
\begin{proof}
  Note that by the choice of the metric on $G$, for any
  $k$-chain $\sigma$ and any $0\le t\le 1$,
  $$\mass s_t(\sigma) \le t^k\mass\sigma$$
  and that for any horizontal $k$-chain $\sigma$ (in particular, when
  $\sigma=P_{\phi(\tau)}(s_{2^{-i}}(\alpha))$) and any $t>0$,
  $$\mass s_t(\sigma)=t^k\mass\sigma.$$
  By Theorem~\ref{thm:fedflem}, there is a $c$ such that
  $$\mass P_\tau(\sigma)\le c \mass\sigma$$
  for all $k$-chains $\sigma$.  

  If we let $c_P:=\Lip(f^{-1})^k \Lip(\phi)^{k}c$, then 
  $$\mass P_i(\alpha) \le 2^{ik}\Lip(f^{-1})^k \Lip(\phi)^{k}c
  2^{-ik}\mass \alpha\le c_P \mass \alpha$$
  as desired.
\end{proof}

Next, we construct simplicial chains interpolating between two
different approximations of a cycle.  The basic idea is that if
$(\tau_0,f_0:\tau_0\to G)$ and $(\tau_1,f_1:\tau_1\to G)$ are two
triangulations of $G$, we can connect approximations in $\tau_0$ and
$\tau_1$ using a triangulation of $X\times [0,1]$ which interpolates
between $\tau_0$ and $\tau_1$.

Let $(\eta, g:\eta \to G\times [0,1])$ be a triangulation of
$G\times [0,1]$.  For $i=0,1$, suppose that $g^{-1}(G\times\{i\})$ is
a subcomplex of $\eta$ which is isomorphic to $\tau_i$ under an
isomorphism $\iota_i:\tau_i\cong g^{-1}(G\times\{i\})$ such that
$g\circ \iota_i=f_i$.  We say that $\eta$ \emph{restricts to $\tau_i$ on
$G\times\{i\}$}.  Let $\psi:\eta\to G$ be a $(k+1)$-horizontal map and
let $\phi_i:\tau_i\to G$ be defined by $\phi_i=\psi \circ \iota_i$ for
$i=0,1$.  We will construct a horizontal $(k+1)$-chain interpolating
between $P_{\phi_0(\tau_0)}(\alpha)$ and $P_{\phi_1(\tau_1)}(\alpha)$.

\begin{lemma}\label{lem:Rprops}
  There is a $c_R$ depending only on $\psi$ and $\eta$ such that for
  all $i\ge 0$ and for all $k$-cycles $\alpha$, there is a
  $(k+1)$-chain $R_{\psi(\eta)}(\alpha)$ such that
  $$\partial
  R_{\psi(\eta)}(\alpha)=P_{\phi_1(\tau_1)}(\alpha)-P_{\phi_0(\tau_0)}(\alpha)$$
  and
  $$\mass R_{\psi(\eta)}(\alpha) \le  c_R \mass \alpha.$$
  Furthermore, this chain is the image under $\psi_\sharp$ of a
  simplicial chain in $\eta$, so it is horizontal.
\end{lemma}
\begin{proof}
  Let
  $$X(\alpha):=Q_{\tau_1}[(f_1^{-1})_\sharp(\alpha)] +
  g^{-1}_\sharp(\alpha\times [0,1]) -
  Q_{\tau_0}[(f_0^{-1})_\sharp(\alpha)],$$ and note that
  $$\partial
  X(\alpha)=P_{\tau_1}[(f_1^{-1})_\sharp(\alpha)]-P_{\tau_0}[(f_0^{-1})_\sharp(\alpha)].$$
  Define 
  $$R_{\psi(\eta)}(\alpha):=\psi_\sharp(P_\eta(X(\alpha))).$$

  This is the image of a simplicial chain, and 
  $$\partial
  R_{\psi(\eta)}(\alpha)=\psi_\sharp(\partial X(\alpha))=P_{\phi_1(\tau_1)}(\alpha)-P_{\phi_0(\tau_0)}(\alpha),$$
  as desired.  For the bound on the mass of $R_{\psi(\eta)}(\alpha)$,
  note that by Theorem~\ref{thm:fedflem}, there is a $c$ such that
  $$\mass P_\tau(\sigma)\le c \mass\sigma$$
  and
  $$\mass Q_\tau(\sigma)\le c \mass\sigma$$
  for all $k$-chains or $(k+1)$-chains $\sigma$.  Thus
  $$\mass X(\alpha) \le [2 c
  (\Lip g^{-1})^k+(\Lip g^{-1})^{k+1}]\mass \alpha.$$
  If we let 
  $$c_R:=c [2 c
  (\Lip g^{-1})^k+(\Lip g^{-1})^{k+1}]\Lip (\psi)^{k+1},$$
  then 
  $$\mass R_{\psi(\eta)}(\alpha) \le  c_R \mass \alpha.$$
\end{proof}

In this case, we would like to connect $P_0(\alpha)=P_{\phi(\tau)}(\alpha)$ and
$$P_1(\alpha)=s_{2}( P_{\phi(\tau)}(s_{2^{-1}}(\alpha))).$$
Let $(\tau_0,f_0:\tau_0\to G)=(\tau,f)$ and $\phi_0=\phi$, and define
$(\tau_1,f_1:\tau_1\to G)$ by letting $\tau_1=\tau$ and $f_1=s_2\circ
f$.  If we define $\phi_1=s_2\circ \phi$, then we have
$P_1(\alpha)=P_{\phi_1(\tau_1)}(\alpha).$ We will define
$R_0(\alpha)$ as $R_{\psi(\eta)}(\alpha)$ for an appropriate $\psi$
and $\eta$ and obtain $R_i$ by conjugating $R_0$ by $s_{2^i}$.
\begin{lemma}\label{lem:constRi}
  Let $k>0$.  Let $(\tau, f:\tau\to G)$ be a triangulation
  and let $\phi:\tau\to G$ be a $(k+1)$-horizontal map which is a
  bounded distance from $f$.  Define $\tau_0$, $\phi_0$, $\tau_1$ and
  $\phi_1$ as above.

  Let $(\eta, g:\eta\to G\times[0,1])$ be a triangulation of $G\times
  [0,1]$ which restricts to $\tau_i$ on $G\times \{i\}, i=0,1$.  Let
  $\iota_i:\tau_i\cong g^{-1}(G\times \{i\})$ be the implied
  isomorphism.  Let $\psi:\eta\to G$ be a $(k+1)$-horizontal map which
  extends the $\phi_i$ (i.e., $\phi_i=\psi\circ \iota_i, i=0,1$).  If
  we define
  $$R_i(\alpha):=s_{2^i}( R_{\psi(\eta)}(s_{2^{-i}}(\alpha))),$$ 
  then for all $i\ge 0$ and for all $k$-cycles $\alpha$, we have
  $$\partial
  R_i(\alpha)=P_{i+1}(\alpha)-P_{i}(\alpha)$$
  and
  $$\mass R_i(\alpha) \le  c_R 2^i \mass \alpha,$$
  where $c_R$ is the constant from Lemma~\ref{lem:Rprops}
  corresponding to $\psi$ and $\eta$.
\end{lemma}
\begin{proof}
  It follows from Lemma~\ref{lem:Rprops} that
  \begin{align*}
    \partial R_{i}(\alpha) & =s_{2^i}[P_{\phi_1(\tau_1)}(s_{2^{-i}}(\alpha))-P_{\phi_0(\tau_0)}(s_{2^{-i}}(\alpha))]\\
    & =(s_{2^i}\circ P_{\phi_1(\tau_1)}\circ
    s_{2^{-i}})(\alpha)-(s_{2^i}\circ P_{\phi_0(\tau_0)} \circ s_{2^{-i}})(\alpha).
  \end{align*}
  Since  
  $$P_{\phi_1(\tau_1)}(\alpha)=(s_2\circ P_{\phi_0(\tau_0)} \circ s_{2^{-1}})(\alpha),$$
  we have $\partial R_{i}(\alpha) =P_{i+1}(\alpha)-P_i(\alpha)$ as
  desired.

  Next we bound the mass of $R_i(\alpha)$.  Recall that $\mass R_0(\alpha)\le c_R \mass
  \alpha$ for all $\alpha$.  Then
  $$\mass R_i(\alpha)=\mass s_{2^i}(R_i(s_{2^{-i}}(\alpha))).$$
  Since $R_i(s_{2^{-i}}(\alpha))$ is a horizontal $(k+1)$-cycle,
  we have 
  $$\mass R_i(\alpha)\le 2^{(k+1)i} c_R 2^{-ik}\mass \alpha=c_R 2^i\mass
  \alpha,$$
  as desired.
\end{proof}

When such an $R_i$ exists, we can use it to prove filling
inequalities.
\begin{thm}\label{thm:triToFilling}
  Let $k$, $(\tau,f)$, $\phi$, $(\eta,g)$, and $\psi$ satisfy the
  hypotheses of Lemma~\ref{lem:constRi}.  Then
  $$\FV^{k+1}_G(V)\preceq V^{\frac{k+1}{k}}.$$
\end{thm}
\begin{proof}
  It suffices to show that there is a $c$ such that if $\alpha$ is a
  $k$-cycle with sufficiently large volume, then there is a chain
  $\beta$ such that $\partial \beta=\alpha$ and
  $$\mass \beta\le c (\mass \alpha)^{\frac{k+1}{k}}.$$

  First, we claim that when $i$ is large, then $P_{i}(\alpha)=0$.  Let
  $$c'=c_\tau \Lip(f^{-1})^k(\vol \Delta^k)^{-1},$$
  where $c_\tau$ is the constant from Theorem~\ref{thm:fedflem} and
  $\Delta^k$ is the standard euclidean $k$-simplex.
  Let $i_0$ be the integer such that
  $$2^{(i_0-1) k} \le c' \mass \alpha < 2^{i_0 k},$$
  and suppose that $i\ge i_0$.

  Let $X=(P_\tau\circ f^{-1}_\sharp\circ s_{2^{-i}})(\alpha),$ so that
  $P_i(\alpha)=s_{2^{i}}(\phi_\sharp(X))$.  We claim that $X=0$.
  Since $X$ is an integral simplicial $k$-cycle, it suffices to show
  that $\mass X< \vol \Delta^k$.  But, by Theorem~\ref{thm:fedflem},
  we have
  $$\mass X \le c_\tau \Lip(f^{-1})^k 2^{-ki} \mass \alpha<\vol
  \Delta^k,$$ 
  so $X=0$ and thus $P_i(\alpha)=0$.

  We claim that
  $$\beta:=-\bigr(Q_{\phi(\tau)}(\alpha)+\sum_{i=0}^{i_0-1} R_i(\alpha)\bigl)$$
  is a filling of $\alpha$ with mass $\preceq (\mass \alpha)^{(k+1)/k}$.
  First, note that 
  $$\partial\beta=\alpha-P_{i_0}(\alpha)=\alpha.$$ 
  Furthermore, by Lemmas~\ref{lem:constRi} and \ref{lem:Qprops}, we
  have
  \begin{align*}
    \mass \beta & \le c_Q\mass \alpha + \sum_{i=0}^{i_0-1} c_R 2^i \mass
    \alpha \\
    & \le (c_Q + c_R 2^{i_0}) \mass \alpha \\
    & \le (c_Q + 2 c_R (c'\mass \alpha)^{1/k}) \mass \alpha.
  \end{align*}
  If $\mass \alpha$ is sufficiently large, then we have
  $$\mass \beta \le 4 c_R (c'\mass \alpha)^{\frac{k+1}{k}}$$
  as desired.
\end{proof}

So we can get filling inequalities by constructing triangulations and
$k$-horizontal maps.  In the next section, we will describe two ways
of constructing such triangulations.

\subsection{Homotopic filling bounds}
In this section, we give a sketch of how to adapt the above arguments
to produce fillings of spheres by balls rather than fillings of chains
by cycles.  Note that when $d\ge 2$, an upper bound on $\FV^{d+1}$
implies a bound on $\delta^d$; Gromov \cite[App.  2.(A')]{GroFRM} and
White \cite{WhiteMinimal} showed that $\delta^d\preceq \FV^{d+1}$
(with equality when $d\ge 3$).

The main change to produce homotopic fillings rather than homological
fillings is a homotopic version of the Deformation Theorem.  To state
this, we recall the \emph{admissable maps} used by Brady, Bridson,
Forester, and Shankar \cite{BBFS}.  If $M$ is a compact
$k$-dimensional manifold and $\tau$ is a simplicial complex, a map
$\alpha:M\to \tau$ is admissible if the image of $\alpha$ lies in the
$k$-skeleton $\tau^{(k)}$ of $\tau$, and if
$\alpha^{-1}(\tau^{(k)}\setminus \tau^{(k-1)})$ is a disjoint union of
open $k$-balls, each mapped homeomorphically onto a $k$-cell of
$\tau$.  If $(\tau,f:\tau\to N)$ is a triangulation of $N$, we say
that $\alpha:M\to N$ is $\tau$-admissible if and only if $f^{-1}\circ
\alpha$ is an admissible map to $\tau$.  One can then prove the
following:
\begin{thm}[Homotopic Deformation
  Theorem]\label{thm:htpicfedflem}
  Let $(\tau,f:\tau\to X)$ be a triangulation of $X$ and let $M$ be
  a compact $k$-manifold.
  
  There is a constant $c=c(\tau)$ such that if $\alpha:M\to X$ is a
  Lipschitz map and $\alpha|_{\partial M}$ is $\tau$-admissible, then
  there is a $\tau$-admissible map $P'_\tau(\alpha):M\to X$ and a
  Lipschitz homotopy $Q'_\tau(\alpha):M\times [0,1] \to X$ such that
  \begin{enumerate}
  \item $\vol P'_\tau(\alpha)\le c \cdot\vol \alpha$,
  \item $\vol Q'_\tau(\alpha)\le c \cdot\vol \alpha$, and
  \item $Q'_\tau(\alpha)$ is a homotopy between $\alpha$ and
    $P'_{\tau}(\alpha)$ which is constant on $\partial M$.
  \end{enumerate}
\end{thm}

This can be used to prove:
\begin{thm}\label{thm:htpicTriToFilling}
  Let $k$, $(\tau,f)$, $\phi$, $(\eta,g)$, and $\psi$ satisfy the
  hypotheses of Lemma~\ref{lem:constRi}.  Then
  $$\delta^{k}_G(V)\preceq V^{\frac{k+1}{k}}.$$
\end{thm}
\sketch If $\alpha:S^k\to X$, we can use the Homotopic Deformation
Theorem to construct $(\tau,s_{2^i}\circ f)$-admissible maps
$P'_i(\alpha):S^k\to X$ approximating $\alpha$ at different scales and
horizontal homotopies $R'_i(\alpha)$ which connect the
$P'_i(\alpha)$'s.  When $k$, $(\tau,f)$, $\phi$, $(\eta,g)$, and
$\psi$ satisfy the hypotheses of Lemma~\ref{lem:constRi}, these maps
and homotopies satisfy the same bounds on their volumes as their
homological counterparts.

We construct a disc filling $\alpha$ by connecting the $R'_i(\alpha)$
end-to-end for $i=0,\dots,i_0$; this gives us a homotopy between
$P'_0(\alpha)=P'_\tau(\alpha)$ and $P'_{i_0}(\alpha)$.  The map
$P'_0(\alpha)$ is homotopic to $\alpha$ by the homotopy
$Q'_\tau(\alpha)$.  Since $P'_{i_0}(\alpha)$ is $(\tau,s_{2^i}\circ
f)$-admissible, the map $\alpha'=f^{-1}\circ s_{2^{-i}}\circ
P'_{i_0}(\alpha)$ is admissible, but if $i_0$ is sufficiently large,
$\vol \alpha'< \vol \Delta^k$.  This implies that the image of
$\alpha'$ lies entirely in the $(k-1)$-skeleton of $\tau$, so $\vol^k
P'_{i_0}(\alpha)=0$.  There is thus a homotopy contracting
$P'_{i_0}(\alpha)$ to a point which has no $(k+1)$-volume.  Combining
all of these homotopies, we get a homotopy from $\alpha$ to a point
whose volume is $\preceq (\vol \alpha)^{(k+1)/k}$.  This homotopy is
the required disc filling $\alpha$.

\section{Constructing horizontal triangulations}\label{sec:horizTris}
Theorems~\ref{thm:triToFilling} and \ref{thm:htpicTriToFilling} show that when certain triangulations
and $k$-horizontal maps exist, then $G$ satisfies euclidean filling
inequalities.  We will give two families of examples of groups with
such maps and triangulations, one using combinatorial group theory and
one using geometry.  The combinatorial example gives a way to prove
quadratic bounds on the Dehn functions of central products; these
groups include the higher-dimensional Heisenberg groups and have been
studied by Ol'shanskii and Sapir \cite{SapOls}, Allcock
\cite{Allcock}, and Magnani \cite{Magnani}.  The geometric example
shows that jet groups have enough horizontal submanifolds to construct
the required maps and triangulations, and that consequently, their
filling functions satisfy euclidean bounds.

\subsection{Central products}
If $G$ is a two-step nilpotent group with center $Z$, we define the
central product $G\times_Z G$ to be the quotient $G\times G/\sim$,
where $\sim$ is the relation which identifies the centers of the two
copies of $G$ (i.e., the relation $(z,1)\sim (1,z)$ for all $z\in Z$).
Likewise, we define the $n$-fold central product to be the quotient 
$$G^{\times_Z n}=G^n/\sim,$$
where $\sim$ is the relation 
$$(z,1,1,\dots)\sim (1,z,1,\dots)\sim(1,1,z,\dots)\sim\dots\text{ for
  all $z\in Z$}$$
which identifies the centers of all the copies of $Z$.  These
groups include the higher-dimensional Heisenberg groups, which are
given by $H_{2n+1}=(H_3)^{\times_Z n}$. 
Central products also appear as cusp groups in
lattices in rank-1 symmetric spaces.  

In this section, we will prove Dehn function bounds for such groups.
First, we will show that central products of free nilpotent groups
have quadratic Dehn functions.  A central product of any two-step
nilpotent group is a quotient of a central product of a free nilpotent
group.  We will show that any such quotient has a Dehn functions which
grows at most as fast as $n^2\log n$ and that many such quotients have
quadratic Dehn functions.

Consider the free two-step nilpotent group on $r$ generators, which is
given by the presentation:
$$\Lambda_r=\langle g_1,\dots, g_r \mid [g_i,[g_j,g_k]]=0 \text{ for
  $1\le i,j,k\le r$.}\rangle.$$ Its abelianization is $\Z^r$,
generated by the $g_i$, and its center is isomorphic to
$\Z^{\binom{r}{2}}$, generated by elements of the form $[g_i,g_j],
1\le i<j\le r$.  Define $\Lambda_{r,n}=(\Lambda_r)^{\times_Z n}$.  
We claim:
\begin{prop}\label{prop:lambdarnQuad}
  $\Lambda_{r,n}$ has a quadratic Dehn function.
\end{prop}
This proposition was first stated without proof by Ol'shanskii and Sapir \cite{SapOls}.  It
will follow from applying Theorem~\ref{thm:htpicTriToFilling} with
$k=1$.  We will show that $\Lambda_{r,n}$ satisfies the conditions of
the theorem by finding a presentation of $\Lambda_{r,n}$ and then
using that presentation to construct triangulations and horizontal
maps.

These groups are lattices in nilpotent Lie groups:  $\Lambda_r$ is
a lattice in the free 2-step nilpotent Lie group of rank $r$, which we
call $F_r$, and $\Lambda_{r,n}$ is a lattice in
$F_{r,n}:=(F_r)^{\times_Z n}$.  Let $\frak{f}_r$ be the Lie algebra of
$F_r$ and let $v_i=\log g_i \in \frak{f}_r$.  If we define generators
of $\frak{f}_r$ by $v_i=\log g_i\in \frak{f}_r$, then
$\frak{f}_r$ has a grading
$$\frak{f}_r=V_1\oplus V_2=\langle v_1,\dots, v_r\rangle \oplus \langle [v_i,v_j]
\text{ for all $1\le i<j\le r$}\rangle.$$ If $g_{ij}$ is the $i$th
generator of the $j$th factor in the central product $F_{r,n}$, and
$\frak{f}_{r,n}$ is its Lie algebra, we can likewise define a grading
$$\frak{f}_{r,n}=V_1^n\oplus V_2=\langle v_{11},\dots, v_{rn}\rangle \oplus
\langle [v_{i1},v_{j1}] \text{ for all $1\le i<j\le r$}\rangle$$
where $v_{ij}=\log g_{ij}$.

For ease of notation, we will start by considering $\Lambda_{r,2}$.  Let
$g_1,\dots, g_r\in \Lambda_r\times \Lambda_r$ be the generators of the
first factor and let $h_1,\dots, h_r\in \Lambda_r\times \Lambda_r$ be
the generators of the second.  Then $\Lambda_{r,2}$ is
given by the presentation:
\begin{align}
  \label{eq:gClass2}\Lambda_{r,2}=\langle g_1,\dots, g_r \mid & [g_i,[g_j,g_k]]=0 \text{ for
    $1\le i,j,k\le r$,}\\
  \label{eq:hClass2}  & [h_i,[h_j,h_k]]=0 \text{ for $1\le i,j,k\le r$},\\
  \label{eq:product}  & [g_i,h_j]=0 \text{ for $1\le i,j\le r$},\\
  \label{eq:centerAmalg}  & [g_i,g_j]=[h_i,h_j]=0 \text{ for  $1\le i,j\le r$}\rangle.
\end{align}
We claim:
\begin{lemma}\label{lem:lambdar2}
  $\Lambda_{r,2}$ can be presented as:
  \begin{align}
\label{eq:product2}    \Lambda_{r,2} =\langle g_1,\dots, g_r \mid & [g_i,h_j]=0 \text{ for $1\le i,j\le r$,}\\
\label{eq:commute}    &[g_ih_j,g_jh_i]=0 \text{ for
      $1\le i,j\le r$}\rangle.
  \end{align}
\end{lemma}
\begin{proof}
  We need to show that the relation $[g_ih_j,g_jh_i]=0$ holds in
  $\Lambda_{r,2}$ for all $i$ and $j$, and we need to show
  that relations \eqref{eq:gClass2}, \eqref{eq:hClass2}, and \eqref{eq:centerAmalg} can
  be deduced from \eqref{eq:product2} and \eqref{eq:commute} (since
  \eqref{eq:product} and \eqref{eq:product2} are the same).

  First, we reduce $[g_ih_j,g_jh_i]$ to the empty word $\varepsilon$ by
  using \eqref{eq:gClass2}--\eqref{eq:centerAmalg}:
  \begin{align*}
    [g_ih_j,g_jh_i] &\to [g_i,g_j][h_j,h_i] & & \text{by
      \eqref{eq:product}} \\
    &\to \varepsilon & &\text{by
      \eqref{eq:centerAmalg}}.
  \end{align*}
  In the first step, we shuffled all the
  $g_i$'s and $g_j$'s in $[g_ih_j,g_jh_i]$ to the beginning using the
  fact that $g$'s and $h$'s commute.

  Now we deduce \eqref{eq:gClass2}, \eqref{eq:hClass2}, and
  \eqref{eq:centerAmalg} from \eqref{eq:product2} and
  \eqref{eq:commute}.  First, note that \eqref{eq:centerAmalg} follows
  from \eqref{eq:product2} and \eqref{eq:commute} by the reverse of
  the argument above:
  $$[g_i,g_j][h_j,h_i] \to [g_ih_j,g_jh_i] \to \varepsilon$$
  Next, we can reduce \eqref{eq:gClass2} as follows:
  \begin{align*}
    [g_i,[g_j,g_k]] &\to [g_i,[h_j,h_k]] & &\text{by
      \eqref{eq:centerAmalg}} \\
    &\to \varepsilon & & \text{by \eqref{eq:product2}}.
  \end{align*}
  Here, we reduced $[g_i,[h_j,h_k]]$ to the trivial word by using the
  fact that $g_i$ commutes with each letter of $[h_j,h_k]$.  The same
  procedure can be used to deduce \eqref{eq:hClass2}.
\end{proof}

Larger central products have a similar presentation; here, $g_{ij}$
represents the $i$th generator of the $j$th factor of $\Lambda_r^n$:
\begin{lemma}
  $\Lambda_{r,n}$ can be presented as:
  \begin{align}
\label{eq:lambdarnPres}    \Lambda_{r,n} =\langle g_{ij}, i=1,\dots, r, j=1,\dots,
    n\mid & [g_{ij},g_{kl}]=0 \text{ for all $j\ne l$,}\\
\notag    &[g_{ij} g_{kl},g_{kj} g_{il}]=0 \text{ for
      all $j\ne l$}\rangle.
  \end{align}
\end{lemma}
We omit the proof, which follows the outline of Lemma~\ref{lem:lambdar2}.

Recall that if $\Gamma$ is a group with a finite presentation
$$\Gamma=\langle x_1,\dots, x_d\mid r_1,\dots, r_s\rangle,$$
then its {\em Cayley complex} $X_\Gamma$ is a simply-connected
2-complex on which $\Gamma$ acts geometrically (that is, properly
discontinuously, cocompactly, and by isometries).  The 1-skeleton of
$X_\Gamma$ is the Cayley graph of $\Gamma$ with respect to the $x_i$.
At each vertex of the Cayley graph, there is a loop corresponding to
each relator $r_i$, and we obtain $X_\Gamma$ by gluing a 2-cell to
each such loop.  Since we started with a presentation of $\Gamma$,
this procedure results in a simply-connected complex.

The advantage of the presentation in \eqref{eq:lambdarnPres} is that
there is a horizontal map from its Cayley complex
$X:=X_{\Lambda_{r,n}}$ to $F_{r,n}$, which we will denote $h:X\to
F_{r,n}$.  Vertices of $X$ correspond to elements of $\Lambda_{r,n}$,
so we map each vertex to the corresponding element of $F_{r,n}$.  Each
edge $e$ of $X$ connects $g$ and $gg_{ij}^{\pm1}$ for some $g\in
\Lambda_{r,n}$ and some $i$ and $j$.  Since $g_{ij}=\exp v_{ij}$ for
all $i,j$, the points $g$ and $gg_{ij}^{\pm1}$ can be connected by a
horizontal segment of the form $t\mapsto g \exp \pm tv_{ij}, 0\le t\le
1$; we define $h(e)$ to be this segment.  Then if $w=w_1\dots w_p$ is
a word, it corresponds to an edge path in $X$ and its image under $h$
is a horizontal curve $\gamma_w$ in $G$ which connects the points
$0,w_1,w_1w_2,\dots, w_1\dots w_p$.  Each relation in
\eqref{eq:lambdarnPres} then corresponds to a horizontal closed curve
$\gamma_{w}$.  We complete the definition of $h$ by filling these
curves with horizontal discs.

\begin{figure}
  \begin{center}
    \def\svgwidth{.5\textwidth}
    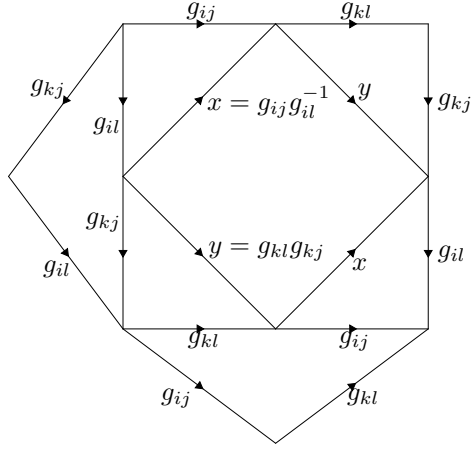
  \end{center}
  \caption{A disc filling $\gamma_w$ for $w=[g_{ij} g_{kl},g_{kj} g_{il}]$\label{fig:extendSq}}
\end{figure}

First, we consider $w=[g_{ij},g_{kl}]$.  The curve $\gamma_w$ lies in the
2-parameter subgroup $\exp\langle v_{ij},v_{kl}\rangle$, and we can
fill it with a disc of the form $\exp(s v_{ij}+t v_{kl}), 0\le s,t\le
1$.  Next, we consider $w=[g_{ij} g_{kl},g_{kj} g_{il}]$.  Here, the
disc is a little more complicated and is shown in
Figure~\ref{fig:extendSq}.
The disc is made up of three quadrilaterals and four triangles.  Each
edge in the disc is a segment of a translate of a horizontal
1-parameter subgroup, and each face lies in a translate of a
2-dimensional horizontal subgroup of $F_{r,n}$.
The quadrilateral on the left lies in a translate of the subgroup with
Lie algebra $\langle v_{il},v_{kj}\rangle$, the one on the bottom lies
in a translate of $\exp \langle
v_{kl},v_{ij}\rangle$, and the one in the center lies in a translate
of $\exp \langle
v_{ij}-v_{il},v_{kl}+v_{kj}\rangle$.  Likewise, each triangle lies in
a translate of either $\exp \langle v_{ij},v_{il}\rangle$ or $\exp
\langle v_{kj},v_{kl}\rangle$.  All of these subgroups are horizontal
submanifolds of $F_{r,n}$, so the disc is horizontal.

As a consequence, we conclude the following:
\begin{lemma}\label{lem:lambdarnFill}
  If $w$ is a word in $\Lambda_{r,n}$ which represents the identity,
  then there is a horizontal disc filling $\gamma_w$ in $F_{r,n}$.
\end{lemma}
\begin{proof}
  Since $w$ represents the identity, it corresponds to a closed edge
  path in $X$, and since $X$ is simply connected, this path can be
  filled by a disc in $X$.  The image of this disc in $F_{r,n}$ is a
  disc filling $\gamma_w$.
\end{proof}

We can use this lemma to construct triangulations and maps satisfying
the hypotheses of Theorem~\ref{thm:htpicTriToFilling} with $k=1$.  To apply
the theorem to a group $G$, we need a triangulation $\eta$ of $G\times
[0,1]$ which restricts to triangulations $\tau_0$ and $\tau_1$ on
$G\times \{0\}$ and $G\times\{1\}$ respectively, where $\tau_1$ is a
scaling of $\tau_0$.  In the following lemma, we show that if $G$
contains a lattice $\Gamma$ such that $\Gamma\subset s_2(\Gamma)$,
then $\tau_0$ can be chosen to be $\Gamma$-equivariant and $\eta$ can
be chosen to be $s_2(\Gamma)$-equivariant.

\begin{lemma}\label{lem:trisExist}
  Let $\Gamma$ be a lattice in a Carnot group $G$ such that
  $s_2(\Gamma)\subset \Gamma$ and let $(\tau, f:\tau \to G)$ be a
  $\Gamma$-equivariant triangulation.  Let $(\tau_i, f_i:\tau_i\to
  G)$, $i=0,1$ be triangulations such that $\tau_0=\tau_1=\tau$,
  $f_0=f$, and $f_1=s_2\circ f_0$.  Then there is a
  $s_2(\Gamma)$-equivariant triangulation $(\eta, g:\eta\to
  G\times[0,1])$ of $G\times[0,1]$ which restricts to $\tau_i$ on
  $G\times\{i\}$ for $i=0,1$.
\end{lemma}
\begin{proof}
  First, we construct $s_2(\Gamma)$-actions on $\tau_0$ and $\tau_1$
  which make them equivariant triangulations.  Let $\rho(g):\tau\to
  \tau, g\in \Gamma$ be the action of $\Gamma$ on $\tau$.  Define
  $\rho_0(g)=\rho(g)$ and $\rho_1(g)=\rho(s_{1/2}(g))$ for all $g\in
  s_2(\Gamma)$.  Then $f_0$ and $f_1$ are equivariant with respect to
  $\rho_0$ and $\rho_1$ respectively.

  In particular, $\tau_0$ and $\tau_1$ descend to triangulations of
  $M:=s_2(\Gamma)\backslash G$.  Since this is a smooth manifold, we
  can construct a triangulation $(\bar{\eta},\bar{g}:\bar{\eta}\to
  M\times [0,1])$ which restricts to $s_2(\Gamma)\backslash\tau_i$ on
  $M\times\{i\}$.  This lifts to the required triangulation of $G\times[0,1]$.
\end{proof}

Since $g_{ij}=\exp v_{ij}$ and $v_{ij}$ is in the first layer of the
grading of $\frak{f}_{r,n}$, we have $s_2(g_{ij})=\exp
2v_{ij}=g_{ij}^2$ and thus $s_2(\Lambda_{r,n})\subset \Lambda_{r,n}$.
The lemma thus applies to $\Lambda_{r,n}$, and we can let $\tau$ and
$\eta$ be the triangulations of $F_{r,n}$ and $F_{r,n}\times [0,1]$
which are given by the lemma.

Now we can define $\phi$.  For each vertex $v$ of $\tau$, let
$\phi(v)$ be an element of $\Lambda_{r,n}$; we can choose these
elements in a $\Lambda_{r,n}$-equivariant way.  If $e$ is an edge of
$\tau$ with vertices $x$ and $y$, we can choose a word
$w=w(\phi(x)^{-1}\phi(y))$ in $\Lambda_{r,n}$ which represents
$\phi(x)^{-1}\phi(y)$.  Then $\gamma_w$ is a curve connecting $e$ to
$\phi(x)^{-1}\phi(y)$, and we can define $\phi(e)$ to be the
translation $\phi(e)=\phi(x)\cdot \gamma_{w}$.  Since $w$ depends only
on $\phi(x)^{-1}\phi(y)$, this definition is also equivariant.  The
map $\phi$ then sends the boundary of each $2$-cell $\Delta$ of $\tau$
to a curve $g\cdot \gamma_{w'(\Delta)}$, where $w'(\Delta)$ is a word
representing the identity.  By Lemma~\ref{lem:lambdarnFill}, we can
extend $\phi$ to a horizontal map on $\Delta$, and furthermore, we can
do this in an equivariant way.  Doing this for every cell in the
2-skeleton of $\tau$ gives us an equivariant horizontal map on
$\tau^{(2)}$, and we can extend it to an equivariant $2$-horizontal
map on all of $\tau$.

Next, we define $\psi$ in the same fashion.  This time, because of the
relationship between $\tau$ and $\eta$, many of our choices are
already made for us.  Let $\tau_i$ and $f_i$, $i=0,1$ be as in
Lemma~\ref{lem:trisExist} and let $\iota_i:\tau_i\to \eta$, $i=0,1$ be
the inclusions of the $\tau_i$ into $\eta$.  Let
$$\psi|_{\iota_0(\tau_0)}=\phi\circ \iota_0^{-1}$$
and
$$\psi|_{\iota_1(\tau_1)}=s_2\circ \phi\circ \iota_1^{-1}.$$
These definitions are $s_2(\Lambda_{r,n})$-equivariant and if $v$ is a
vertex of $\iota_i(\tau_i),i=0,1$, then $\phi(v)\in \Lambda_{r,n}$.
We extend
$\psi$ to the rest of the 1-skeleton of $\eta$ in the same way as
before.  For each vertex $v$ of $\eta$ which is not in $\tau_0$ or
$\tau_1$, we choose an element $\psi(v)\in \Lambda_{r,n}$.  For each
edge $e=(x,y)$ of $\eta$ which is not in $\tau_0$ or $\tau_1$, we let
$$\psi(e)=\psi(x)\cdot\gamma_{w(\psi(x)^{-1}\psi(y))}.$$
Then if $e$ is an edge of $\eta$, then $\psi(e)$ is a
reparameterization of $\gamma_w$ for some word $w$; this is true by
construction if $e$ is not in $\tau_1$, and if $e$ is an edge of
$\tau_1$, then $\psi(e)$ is a curve of the form $s_2\circ \gamma_{w}$
for some $w=w_1\dots w_p$.  Since each of the curves making up
$\gamma_w$ is a segment of a horizontal one-parameter subgroup,
$s_2\circ \gamma_{w}$ is a reparameterization of $\gamma_{z}$, where
$z=w_1^2\dots w_p^2$.  Thus $\psi$ sends the boundary of each $2$-cell
$\Delta$ of $\eta$ to a reparameterization of a curve $g\cdot
\gamma_{w'(\Delta)}$, where $w'(\Delta)$ is a word representing the
identity.  Using Lemma~\ref{lem:lambdarnFill}, we can extend $\psi$ to
a horizontal map on $\Delta$.  Finally, we extend $\psi$ to a
$2$-horizontal map on all of $\eta$.  All of the above extensions can
be made equivariantly and Lipschitz-ly, so that $\phi$ and $\psi$
satisfy the conditions of Lemma~\ref{lem:constRi} for $k=1$.  By
Theorem~\ref{thm:htpicTriToFilling}, we
conclude that $F_{r,n}$ and $\Lambda_{r,n}$ both satisfy quadratic
filling inequalities when $n\ge 2$; this proves
Prop.~\ref{prop:lambdarnQuad}. 

As a corollary, we see that the higher-dimensional Heisenberg groups,
$\Lambda_{2,n}, n\ge 2$, have quadratic Dehn functions; this result
was first proved by Gromov \cite{GroAII} and subsequently using
different methods by Ol'shanskii and Sapir \cite{SapOls} and Allcock
\cite{Allcock}.

Other central products of nilpotent groups can be written as quotients
of the $\Lambda_{r,n}$, which lets us prove bounds on their Dehn
functions as well.  The following two propositions concern quotients
of two-step nilpotent groups with quadratic Dehn functions by central
subgroups.  They show that any such quotient has a Dehn function at
most $\ell^2\log\ell$, and that a quotient by a commutator has
quadratic Dehn function.  These results complement a theorem of Wenger
\cite{WengerQuad} that states that quotients of nilpotent Lie groups
by central subgroups that are not generated by commutators cannot have
quadratic Dehn functions.  Indeed, combining Wenger's theorem and
Cor.~\ref{cor:quadQuot} gives a criterion for deciding whether a
quotient of a two-step group with quadratic Dehn function also has
quadratic Dehn function.

In particular, the central products of the quaternionic and octonionic
Heisenberg groups satisfy the hypotheses of Cor.~\ref{cor:quadQuot}.
These groups appear as cusp groups in finite-volume rank-1 locally
symmetric spaces, and Gromov \cite{GroAII} showed that they have
quadratic Dehn functions using microflexibility.  We will use
Cor.~\ref{cor:quadQuot} to give an alternate combinatorial proof of
this fact.

\begin{prop}\label{prop:quadQuot}
  Let $G$ be a 2-step nilpotent Lie group with quadratic Dehn function
  and let $\frak{g}$ be its Lie algebra.  If $w_1,w_2\in \frak{g}$, let
  $G'$ be the nilpotent Lie group with Lie algebra $\frak{g}/\langle
  [w_1,w_2]\rangle$.  Then $G'$ also has quadratic Dehn function.
\end{prop}

\begin{cor}\label{cor:quadQuot}
  Let $G$ be a 2-step nilpotent Lie group with quadratic Dehn function
  and let $\frak{g}$ be its Lie algebra.  If $W$ is generated by
  elements of the form $[w_1,w_2]$ with $w_1,w_2\in \frak{g}$, let $G'$
  be the nilpotent Lie group with Lie algebra $\frak{g}/W$.  Then
  $G'$ has quadratic Dehn function.
\end{cor}

\begin{prop}\label{prop:logQuadQuot}
  Let $G$ be a nilpotent Lie group with quadratic Dehn function and
  let $\frak{g}$ be its Lie algebra.  If $\frak{g}$ has a grading
  $\frak{g}=V_1\oplus V_2$ and if $W\subset V_2$, let $G'$ be the
  nilpotent Lie group with Lie algebra $\frak{g}/W$.  Then 
  $$\delta_{G'}(\ell)\preceq \ell^2\log \ell.$$
\end{prop}

\begin{cor}\label{cor:logQuadQuot}
  Let $G$ be a 2-step nilpotent Lie group and let $n\ge 2$.  If
  $G'=G^{\times_Z n}$, then 
  $$\delta_{G'}(\ell)\preceq \ell^2\log \ell.$$
\end{cor}
\begin{proof}
  If the abelianization of $G$ has rank $d$, we can write $G=F_d/R$
  for some $R\subset [F_d,F_d]$ and apply
  Prop.~\ref{prop:logQuadQuot} to the quotient $G'=F_d^{\times_Z
    n}/R$.
\end{proof}
This corollary was first stated without proof by Ol'shanskii and Sapir
\cite{SapOls}.

The proofs of both of the propositions use the following construction,
which lifts curves in $G'$ to curves in $G$.  Let $G$, $G'$, and $W$
be as in Prop.~\ref{prop:logQuadQuot}, and let $p:G\to G'$ be the
quotient map.  Let $\frak{g}$ and $\frak{g}'$ be the Lie algebras of
$G$ and $G'$ respectively, so that $\frak{g}=V_1\oplus V_2$ and
$\frak{g}'=V_1\oplus V_2/W$.  We distinguish brackets in $\frak{g}$
and $\frak{g}'$ by writing $[\cdot,\cdot]$ or $[\cdot,\cdot]'$
respectively.

Choose metrics on $V_1$ and $V_2$.  These induce left-invariant
metrics on $\frak{g}$ and $\frak{g}'$.  If $W^\perp\subset \frak{g}$
is the orthogonal complement of $W$, then $W^\perp$ is naturally
isometric to $\frak{g}'=\frak{g}/W$, and we can define $i:\frak{g}'\to \frak{g}$
to be that inclusion.  Define $p_W:V_2\to W$ to be the orthogonal
projection to $W$.

For each Lipschitz curve $\gamma:[0,t]\to G'$, we will define a lift
$\tilde\gamma:[0,t]\to G$ in such a way that if $\gamma$ is $C^1$ and
$\dot{\gamma}(t)\in \frak{g}'$ is its tangent vector, we will have
$\dot\tilde\gamma=i(\dot\gamma)$.  Let $\gamma:[0,t]\to G'$ be a
Lipschitz curve and let $\gamma_1:[0,t]\to V_1$ and $\gamma_2:[0,t]\to
V_2/W$ be such that
$$\gamma(t)=\exp(\gamma_1(t)+\gamma_2(t))$$
for all $t$.  If $\gamma$ is $C^1$, then 
$$\dot{\gamma}(t)=\dot{\gamma}_1(t)+\dot{\gamma}_2(t)-\frac{1}{2}[\gamma_1(t),\dot{\gamma}_1(t)]'.$$
Define
$$\tilde{\gamma}(t):=\exp\biggl[\gamma_1(t)+i\bigl(\gamma_2(t)\bigr)+\frac{1}{2}\int_0^tp_W([\gamma_1(x),\dot{\gamma_1}(x)])\;dx\biggr].$$
If $\gamma$ is $C^1$, we have
\begin{align*}
  \dot{\tilde{\gamma}}(t)& =\dot{\gamma}_1(t)+i\bigl(\dot \gamma_2(t)\bigr)+
  \frac{1}{2}\bigl(p_W([\gamma_1(t),\dot{\gamma_1}(t)])-[\gamma_1(t),\dot{\gamma}_1(t)]\bigr)\\
  & =\dot{\gamma}_1(t)+i\bigl(\dot \gamma_2(t)\bigr)-
  \frac{1}{2}i\bigl([\gamma_1(t),\dot{\gamma}_1(t)]'\bigr)=i(\dot \gamma(t))
\end{align*}
The following properties are easy to check:
\begin{lemma}
  If $\gamma:[0,t]\to G'$ is a Lipschitz curve, then 
  \begin{itemize}
  \item $p\circ \tilde{\gamma}=\gamma$
  \item $\ell( \tilde{\gamma})=\ell(\gamma)$
  \item If $\gamma:[0,t]\to G'$ is a closed curve such that
    $\gamma(0)=\gamma(t)=0$, then $\tilde{\gamma}$ is a curve
    connecting $0$ and $\exp w$ for some $w\in W$.  There is a $c$
    depending only on $G$, $W$, and the chosen metrics such that
    $\|w\|\le c\ell(\gamma)^2$.
  \end{itemize}
\end{lemma}

\begin{proof}[Proof of Prop.~\ref{prop:quadQuot}]
  Since $\frak{g}$ is a 2-step nilpotent Lie algebra, it has a grading
  $\frak{g}=V_1\oplus V_2$.  Without loss of generality, we may assume
  that $w_1, w_2\in V_1$; if not, we can replace them by their
  projections to $V_1.$ Note that $w_1$ and $w_2$ generate a
  2-dimensional subalgebra of $\frak{g}/\langle [w_1,w_2]\rangle$.

  Let $\gamma:[0,1]\to G'$ be a closed curve such that $\gamma(0)=0$.
  We will construct a filling of $\gamma$.  Note that $\tilde{\gamma}$
  is a curve connecting $0$ and $\exp \alpha[w_1,w_2]$ for some
  $\alpha$.  If $\alpha<0$, we can reverse $\gamma$ to get a curve
  such that $\alpha>0$, so we may assume that $\alpha\ge 0$.

  Define $q:[0,4]\to G'$ to be the boundary of a square of area
  $\alpha$ in the subgroup generated by $w_1$ and $w_2$:
  $$q(t):=\exp \begin{cases}
    \sqrt{\alpha}tw_1 & 0\le t\le 1 \\
    \sqrt{\alpha}(w_1+(t-1)w_2) & 1< t\le 2 \\
    \sqrt{\alpha}(w_1+w_2-(t-2)w_1)& 2< t\le 3 \\
    \sqrt{\alpha}(w_2-(t-3)w_2)& 3< t\le 4 
  \end{cases}$$ Then $\tilde{q}:[0,4]\to G$ is a curve of length
  $4\sqrt{\alpha}=O(\ell(\gamma))$ which connects $0$ to $\exp
  \alpha [w_1,w_2]$.

  If $\gamma_1$ and $\gamma_2$ are curves in $G$ with the same endpoints,
  then the union of $\gamma_1$ and $\gamma_2$ is a closed curve.  This
  has a filling by a disc of area $\preceq
  \delta_G(\ell(\gamma_1)+\ell(\gamma_2))$, and this disc can be
  reparameterized to be a homotopy between $\gamma_1$ and $\gamma_2$
  which leaves their endpoints fixed.

  Since $\tilde{\gamma}$ and $\tilde{q}$ have the same endpoints,
  there is a homotopy $h$ from $\tilde{\gamma}$ to $\tilde{q}$ which
  leaves the endpoints fixed and which has area
  $$\area h\le \delta_G(\ell(
  \tilde{q})+\ell(\tilde{\gamma}))=O(\ell(\gamma)^2).$$ This projects
  to a homotopy $h'$ between $\gamma$ and $q$.  Since $[w_1,w_2]=0$ in
  $G'$, $q$ is the boundary of a square in $G'$ of area $\alpha$ in
  the $w_1w_2$-plane.  We can fill $\gamma$ by gluing this square to
  $h'$, so
  $$\delta_{G'}(\gamma)\le \area h'+\alpha =O(\ell(\gamma)^2)$$
  as desired.
\end{proof}

\begin{proof}[Proof of Prop.~\ref{prop:logQuadQuot}]
  If $\gamma_1:[0,1]\to G$ and $\gamma_2:[0,1]\to G$ are two curves in
  $G$ such that $\gamma_1(0)=\gamma_2(0)=0$, define
  $\gamma_1*\gamma_2:[0,1]\to G$
  to be the curve obtained by concatenating $\gamma_1$ (with domain
  rescaled to $[0,1/2]$) and
  $\gamma_1(1)\cdot \gamma_2$ (with domain $[1/2,1]$); this curve connects 0 to
  $\gamma_1(1)\gamma_2(1)$.  Define $\gamma_1^{*n}$ to be the
  concatenation of $n$ copies of $\gamma_1$
  parameterized so that the $i$th copy has domain $[(i-1)/n,i/n]$.

  Let $\gamma:[0,1]\to G'$ be a closed curve such that $\gamma(0)=0$.
  Without loss of generality, we may assume $\ell(\gamma)\ge 1$.  Then
  $\tilde{\gamma}(1)$ is an element of $\exp W$; let $w=\log
  \tilde{\gamma}(1)$.  For any $i\ge 0$, let
  $\tilde{\gamma}_i:=(s_{2^{-i}}\circ \tilde{\gamma})^{*4^i}$; this too
  is a curve connecting $0$ and $\exp w$, and we will construct a
  homotopy from $\tilde{\gamma}$ to $\gamma_k$ which passes
  through $\gamma_1,\dots, \gamma_{k-1}$.
  
  Since $G$ has a quadratic Dehn function, there is a homotopy $h$
  between $\gamma_0$ and $\gamma_1$ which fixes the endpoints of the
  curves.  Scalings of this homotopy can be combined to get homotopies
  from $\gamma_i$ to $\gamma_{i+1}$.  Indeed, $s_{2^{-i}}\circ h$ is a
  reparameterization of a homotopy from 
  $$s_{2^{-i}}\circ \gamma_0=\gamma_i|_{[0,4^{-i}]}$$
  to 
  $$s_{2^{-i}}\circ \gamma_1=\gamma_{i+1}|_{[0,4^{-i}]}$$
  and we can build a homotopy $h_i$ from $\gamma_i$ to
  $\gamma_{i+1}$ out of $4^i$ copies of $s_{2^{-i}}\circ h$.  This
  homotopy has area at most $4^i 4^{-i}\area h=\area h$.

  Let $i_0\ge 0$ be an integer such that $2^{i_0}\le \ell(\gamma)<
  2^{i_0+1}$ and consider the homotopy $\bar{h}$ from $\gamma_0$ to
  $\gamma_{i_0+1}$ obtained by concatenating $h_0,\dots, h_{i_0-1}$.
  Let $p:G\to G'$ be the quotient map and consider $p\circ \bar{h}$.
  This is a homotopy from $\gamma$ to $(s_{2^{-i_0}}\circ
  \gamma)^{*4^{i_0}}$ with area $\le i_0\area h=O((\ell(\gamma))^2\log
  \ell(\gamma))$.  The curve $s_{2^{-i_0}}\circ \gamma$ is a closed
  curve with  length
  $\le 1$, so it has filling area $\le \delta_{G'}(1)$ and the
  concatenation $(s_{2^{-i_0}}\circ \gamma)^{*4^{i_0}}$ has filling
  area $\le 4^{i_0} \delta_{G'}(1)=O((\ell(\gamma))^2)$.  Combining
  $p\circ \bar{h}$ with a filling of $(s_{2^{-i_0}}\circ
  \gamma)^{*4^{i_0}}$ results in a filling of $\gamma$ with area
  $O((\ell(\gamma))^2\log \ell(\gamma))$.
\end{proof}

Two families of groups that satisfy the conditions of
Cor.~\ref{cor:quadQuot} are the higher-dimensional quaternionic and
octonionic Heisenberg groups; these groups appear as cusp groups in
finite-volume quaternion-hyperbolic and Cayley-hyperbolic manifolds.
These groups can be defined as follows: let $\K=\C,\Quat,\Octon$ be
the complex numbers, quaternions, or octonions and let $\K^*$ be the
elements of $\K$ with real part 0 (the pure imaginary elements of
$\K$).  We can then define a graded nilpotent Lie algebra
$\frak{h}_\K=\K\oplus \K^*$.  Its bracket is given by a bilinear form
$\omega:\K\times \K \to \K^*$, which we define as
$\omega(v,w)=\Im(v\bar{w})$, where $\Im(x)$ is the pure imaginary part
of $x$.  Then $\frak{h}_\K$ is a 2-step graded nilpotent Lie algebra,
and we define $H_\K$ to be the corresponding nilpotent Lie group.  We
call $H_\C$ the Heisenberg group, and call $H_\Quat$ and $H_\Octon$
the quaternionic and octonionic Heisenberg groups.  Gromov
\cite{GroAII} showed that $H_\C$ has a cubic Dehn function; Pittet
\cite{PittetCLB} extended this result to $H_\Quat$, but the Dehn
function of $H_\Octon$ is unknown (it was claimed to be cubic in
\cite{PittetCLB}, but the calculation contained a sign error,
corrected in \cite{LPQuad}).

Central products of the (quaternionic or octonionic) Heisenberg group
are called higher-dimensional (quaternionic or octonionic) Heisenberg
groups.  The corresponding Lie algebras have presentations based on
the multiplication table of $\K$.  For instance, $\Quat$ is generated as a
vector space by unit quaternions $1, i,j,k$, with multiplication table
\begin{center}\begin{tabular}{c|cccc}
  &$1$ &$i$ &$j$ &$k$ \\
  \hline
  $1$&$1$ & $i$ & $j$ & $k$ \\
  $i$ & $i$ & $-1$ & $k$ & $-j$ \\
  $j$ & $j$ & $-k$ & $-1$ & $i$ \\
  $k$ & $k$ & $j$ & $-i$ & $-1$.
\end{tabular}
\end{center}
If we let $\frak{f}_4$ be the free two-step Lie algebra of rank 4 and
denote its generators by $1,i,j$, and $k$, we can
write $\frak{h}_\Quat$ as a quotient of $\frak{f}_4$ by identifying
brackets of the generators in pairs:
$$\frak{h}_\Quat=\frak{f}_4/\langle[1,i]-[j,k],[1,j]-[k,i],[1,k]-[i,j]\rangle.$$
If we denote the generators of $\frak{f}_{4,n}$ by
$1_1,i_1,j_1,k_1,\dots, 1_n, i_n, j_n, k_n$, we can write
$\frak{h}_\Quat^{\times_Z n}$ as a quotient of
$\frak{f}_{4,n}$:
\begin{align*}
  \frak{h}_\Quat^{\times_Z
    n}&=\frak{f}_{4,n}/\langle[1_1,i_1]-[j_2,k_2],[1_1,j_1]-[k_2,i_2],[1_1,k_1]-[i_2,j_2]\rangle\\
  &=\frak{f}_{4,n}/\langle[1_1+j_2,
i_1-k_2],[1_1+k_2,j_1-i_2],[1_1+i_2,k_1-j_2]\rangle.
\end{align*}
Since this is a quotient by commutators,
Cor.~\ref{cor:quadQuot} applies, and $H_\Quat^{\times_Z n}$ has a
quadratic Dehn function.  The group $H_\Octon^{\times_Z n}$ has a similar
presentation and thus also has a quadratic Dehn function.

\subsection{Jet groups}

In \cite[4.1.D, 4.4.A--B]{GroCC} Gromov used infinitesimal
invertibility and the $h$-principle to prove a Lipschitz Extension
Theorem for many spaces, including the $k$-jet bundle.  This theorem
can be used to construct maps and triangulations satisfying
Theorem~\ref{thm:triToFilling} for many groups.  In the case of the $k$-jet
bundle, these maps and triangulations can be constructed fairly
explicitly (see also \cite{Thom}), and in this section, we provide an
elementary construction of such maps and triangulations for jet
groups, a family of groups based on jet bundles which includes the
higher-dimensional Heisenberg groups.

This family of groups has also appeared as a family of non-rigid Carnot groups
\cite{War} and a family of quadratically presented Lie
algebras \cite{Chen}.  Warhurst defined the group by putting a group
structure on the $m$-jet bundle $J^m(\R^k)$; we will give a version of
this construction.

The $m$-jet bundle is a generalization of the cotangent bundle and
is often used to describe differential relations.  A smooth map from
$M$ to $\R$ has a gradient which can be considered as a map from $M$
to the cotangent bundle, $TM^*$.  Likewise, its $m$-th order
derivative can be considered as a map to the $i$th symmetric power of
$TM^*$.  The $m$-jet bundle is a sum of these symmetric powers.  We
will mainly consider the $m$-jet bundle of $\R^k$, which we denote 
by $J^m(\R^k)$.  This is a vector bundle with fiber
$$W=\bigoplus_{i=0}^m W_i,$$ where
$W_i:=S^i\R^{k*}$ is the $i$th symmetric power of $\R^{k*}$ and
$W_0:=S^0\R^{k*}=\R$.  A $C^n$ map $f:\R^k\to\R$ corresponds to a
$C^{n-m}$ section, called a \emph{prolongation}, $j^m(f):\R^k\to
J^m(\R^k)$, given by taking derivatives: the projection to $W_0$
corresponds to the original map, the projection to $W_1=\R^{k*}$ is
the gradient of $f$, and so on.  We will often write $j^m_p(f)$ in
place of $j^m(f)(p)$.  In the Carnot structure that we will construct on
$J^m(\R^k)$, prolongations of smooth functions will be horizontal.  

Using a basis of $\R^k$, we can construct a basis of $W$.  Let
$\{e_1,\dots,e_k\}$ be the standard basis of $\R^k$, and let $\{
e_1^*,\dots,e_k^*\}$ be the corresponding basis of
$\R^{k*}$.  If we let
$$y_{(a_1,\dots,a_k)}=\prod_{i=1}^k (e_i^*)^{a_i},$$
then 
$$\left\{y_{(a_1,\dots,a_k)}\middle| \sum a_i=n\right\}$$
is a basis of $W_n$ and 
$$\left\{y_{(a_1,\dots,a_k)}\middle| \sum a_i\le m\right\}$$
is a basis of $W$.  
Using this basis, we can write $J^m(\R^k)$ as a product
$J^m(\R^k)=\R^k\times W$.  

We can now define the group structure.  First, note that for every
$p\in J^m(\R^k)$, there is a unique polynomial in $k$ variables of
degree at most $m$ whose prolongation passes through $p$.  We call
this polynomial $P(p)$, and we use it to construct an action of $\R^k$ on $W$.  Define
$D^m_x(f)\in W$ to be the first $m$ derivatives of
$f$ at $x$, so that 
$$j^m_x(f)=(x,D^m_{x}(f)).$$
If $x\in \R^k$ and $w\in W$, we let $S_x:W\to W$ be the map
$S_x(w)=D^m_x(P((0,w)))$.  This is an action of $\R^k$ on $W$, and we define
$$(p_1,p_2)(q_1,q_2)=(p_1+q_1,S_{q_1}(p_2)+q_2).$$
This makes $J^m(\R^k)$ a semidirect product of $\R^k$ and $W$.
One can check that $\R^k\times \{0\}$ is a subgroup and that
the translate $p\cdot \left(\R^k\times \{0\}\right)$ is the graph of $j^m(P(p))$.

To describe the Lie algebra $\mathfrak{j}_{m,k}$ of $J^m(\R^k)$,
consider the basis
$$\left\{y_{(a_1,\dots,a_k)}\middle| \sum a_i\le m\right\}\cup \{e_1,\dots,e_k\}$$
of $\mathfrak{j}_{m,k}$.  Calculating brackets, we find that
$$\left[e_i,y_{(a_1,\dots,a_k)}\right]=y_{(a_1,\dots,a_i-1,\dots,a_k)}\text{\qquad if
  $a_i>0$}.$$ and all other brackets are zero.  This Lie algebra is
isomorphic to the Lie algebra corresponding to the model
$\mathcal{M}_{m,k}$ defined by Chen \cite{Chen}.  We can give
$\mathfrak{j}_{m,k}$ the grading
$$\mathfrak{j}_{m,k}=(\R^k\oplus W_m)\oplus W_{m-1}\oplus\dots\oplus W_0.$$
Furthermore, since the structure constants of $\mathfrak{j}_{m,k}$ with respect to
the basis $\{x_i,y_{(a_1,\dots,a_k)}\}$ are rational, $\{\exp
x_1,\dots,\exp x_k\}\cup\left\{\exp y_{(a_1,\dots,a_k)}\big| \sum
  a_i=m\right\}$ generate a lattice in $J^m(\R^k).$ Call this lattice
$\Gamma_{m,k}$.  Since the generators are in $\exp (\R^k\oplus W_m)$,
we have $s_2(\Gamma_{m,k})\subset \Gamma_{m,k}$.

Note that the groups $J^1(\R^k)$ are the $(2k+1)$-dimensional Heisenberg
groups and $J^2(\R^2)$ is the class 3 example given in
\cite{YoungSca}; one isomorphism between them takes
$$a,b,c,d,e,f,g,h$$
to
$$y_{(2,0)},x_1,y_{(1,1)},x_2,y_{(0,2)},-y_{(1,0)},y_{(0,1)},-y_{(0,0)}$$
respectively.

Warhurst showed that the left-invariant plane field corresponding to
$\R^k\oplus W_m$ agrees with the standard contact structure on
$J^m(\R^k)$ \cite{War}.  This gives a way to construct
horizontal submanifolds: if $U$ is an open subset of $\R^k$ and
$f:U\to \R$ is smooth, we define $M_f$ to be the image of $j^m(f):U\to
J^m(\R^k)$.  Then $M_f$ is a horizontal
submanifold:
\begin{lemma}\label{horizsubmans}
  If $U$ is an open subset of $\R^k$ and $f:U\to \R$ is smooth, then
  $M_f$ is a smooth horizontal submanifold of $J^m(\R^k)$.
\end{lemma}
\begin{proof}
  $M_f$ is smooth by the definition of $j^m(f)$, so it just remains to
  show that its tangent plane lies in a translate of $V_1=\R^k\oplus
  W_m$.

  First, note that any translate of a prolongation is still a
  prolongation.  Specifically, if $p=(p_1,p_2)\in J^m(\R^k)$, then 
  $p\cdot M_f=M_g$ for $g:U'\to \R$ defined by
  $$g(x)=f(x-p_1)+P(p)(x),$$
  where $U'=U+p_1=\{u+p_1\mid u\in U\}$ and $P(p)(x)$ is the polynomial used
  above
  to define the group structure on $J^m(\R^k)$.  

  Let $x\in U$.  If we take $p=[j^m_x(f)]^{-1}$ in the previous
  construction, then $g$ is a smooth function which vanishes to $m$th
  order at $0$.  The tangent plane to $M_f$ at $j^m_x(f)$ is then a
  translate of the tangent plane to $M_g$ at $0$.  Since the first $m$
  derivatives of $g$ disappear at $0$, the tangent plane to $M_g$ at
  $0$ is contained in $\R^k\oplus W_m$.  Indeed, if we consider the
  derivative of the $m$th derivative of $g$ as a map from $\R^k\to
  W_m$, then the tangent plane to $M_g$ is the graph of this map.
\end{proof}

Now we define a class of horizontal manifolds coming from the
construction of Lemma \ref{horizsubmans}.
\begin{dfn}
  If $U\subset \R^k$ is an open subset and $f:U\to \R$ is a smooth
  map, we say that any submanifold $Y$ of $M_f$ is {\em holonomic}.

  If $X$ is a complex, a Lipschitz map $f:X\to J^m(\R^k)$ is holonomic
  if and only if its image lies in $M_f$ for some smooth $f:U\to \R$.
\end{dfn}
Lemma \ref{horizsubmans} then implies that holonomic submanifolds and
maps are horizontal.

A cycle in a holonomic submanifold is equipped with a smooth function,
and we can use this to construct a holonomic filling.  That is, if
$f:U\to \R$ is a smooth function defined on an open set and if
$\alpha$ is a singular Lipschitz cycle in $M_f$, then we can construct
a holonomic filling of $\alpha$.  The support of $\alpha$ is compact,
so there is a smooth function $\bar{f}:\R^k\to \R$ which agrees with
$f$ on a neighborhood of the support of $\alpha$.  In particular, we
can write $\alpha$ as $J^m(\bar{f})_\sharp(\alpha_0)$ for some cycle
$\alpha_0$ in $\R^k$.  If $\beta_0$ is a chain in $\R^k$ which fills
$\alpha_0$, then $J^m(\bar{f})_\sharp(\beta_0)$ is a chain in
$J^m(\R^k)$ which fills $\alpha$.

In general, a horizontal map is not necessarily even locally
holonomic, and it can be difficult to fill an arbitrary horizontal map
with a horizontal filling.  For our purposes, it suffices to fill
locally holonomic maps.  To work with such maps, we will define
\emph{augmented maps} which are horizontal maps locally equipped with
smooth functions.
\begin{dfn}
  Let $X$ be a simplicial complex.  Let $p_{\R^k}:J^m(\R^k)\to \R^k$
  be the bundle projection.  An \emph{augmented map} from $X$ to
  $J^m(\R^k)$ is a tuple
  $$(\alpha:X\to J^m(\R^k), \{f_\Delta:V_\Delta\to \R\}_{\Delta\subset
    X})$$ which satisfies two conditions.  First, we require that the
  map $\alpha$ is holonomic on each cell.  That is, the image of a
  cell $\Delta$ is contained in $M_{f_\Delta}$.  Second, we require
  that if $\Delta_1$ is a face of $\Delta_2$, then
  $M_{f_{\Delta_1}}\subset M_{f_{\Delta_2}}$.  Note that this imposes
  compatibility conditions on any pair of faces that intersect,
  because if $\Delta$ and $\Delta'$ intersect, then $M_{f_{\Delta\cap
      \Delta'}}\subset M_{f_{\Delta}}$ and $M_{f_{\Delta\cap
      \Delta'}}\subset M_{f_{\Delta'}}$, so $f_{\Delta}$ and
  $f_{\Delta'}$ agree on a neighborhood of $\Delta\cap\Delta'$.
\end{dfn}

If $X$ is a subcomplex of $Y$ and if
$(\alpha,\{f_\Delta\}_{\Delta\subset X})$ is an augmented map on $X$,
we say that $(\beta,\{g_\Delta\}_{\Delta\subset Y})$ \emph{extends}
$(\alpha,\{f_\Delta\})$ if and only if $\beta$ extends $\alpha$ in the
ordinary sense and $M_{g_\Delta}\subset M_{f_\Delta}$ for all simplices $\Delta
\subset X$.  If $\kappa:X\to \R^k$, then we
say that $(\alpha,\{f_\Delta\})$ \emph{covers} $\kappa$ if
$p_{\R^k}\circ \alpha=\kappa$.

\begin{lemma} \label{lemma:simplexextension}
  If $\kappa:\Delta\to \R^k$ is a Lipschitz embedding of a simplex and
  $$(\alpha,\{f_\delta:V_\delta\to
  \R\}_{\delta\subset\bd\Delta})$$ is an augmented map on $\bd \Delta$
  which covers $\kappa|_{\bd\Delta}$, then there is an augmented
  map $(\beta,\{g_\delta:W_\delta\to\R\}_{\delta\subset\Delta})$
  which extends $(\alpha,\{f_\delta\})$ and covers $\kappa$.
\end{lemma}
\begin{proof}
  It suffices to find a smooth function $g:\R^k\to \R$ such that $g$
  agrees with $f_\delta$ on a neighborhood of $\kappa(\delta)$ for
  each face $\delta$.  If we have such a $g$, we can construct the
  required extension by letting $g_\Delta=g$, $\beta=J^m(g)\circ
  \kappa$, and $W_\Delta=\R^k$.  To find $g$, we use a partition of
  unity.

  If $\delta$ is a face of $\Delta$, let
  $$U_\delta:=V_\delta\setminus \bigcup_{\delta'\not\supset \delta}
  \kappa(\delta') \text{\qquad for all $\delta\subsetneq \Delta$}$$
  and let $U_\Delta=\R^k\setminus \kappa(\partial \Delta)$.  For each
  $\delta$, we have $\kappa(\interior \delta)\subset U_\delta$, so the
  $U_\delta$'s form an open cover of $\R^k$.  Let
  $\{\rho_\delta:\R^k\to \R\}_{\delta\subset \Delta}$ be a smooth
  partition of unity subordinate to $\{U_\delta\}$.

  Let $f_\Delta:\R^k\to \R$ be a smooth function and let $g:\R^k\to
  \R$ be $g=\sum_\delta \rho_\delta f_\delta.$ We claim that if
  $\delta$ is a simplex of $\Delta$, then $g$ agrees with $f_\delta$
  on some neighborhood of $\kappa(\delta)$.  Let $x\in
  \kappa(\delta)$.  If $x\in \supp \rho_{\delta'}$, then $x\in
  U_{\delta'}$.  Furthermore, $U_{\delta'}\cap
  \kappa(\delta)=\emptyset$ when $\delta' \not\subset \delta$, so if
  $x\in \supp \rho_{\delta'}$, then $\delta'\subset \delta$.  There is
  thus a neighborhood $U$ of $x$ such that
  $$g|_U=\sum_{\delta'\subset \delta}\rho_{\delta'} f_{\delta'}.$$
  All of the $f_{\delta'}$ agree with $f_\delta$ on a neighborhood of
  $x$, so $g$ agrees with $f_\delta$ on a neighborhood of $x$.  Since
  this is true for any $x\in \kappa(\delta)$, we conclude that
  $f_\delta$ and $g$ agree on some neighborhood of $\kappa(\delta)$.
  We can thus construct the required extension by letting $g_\Delta=g$
  and $\beta=J^m(g)\circ \kappa$.
\end{proof}

This lets us construct maps satisfying the conditions of
Theorem~\ref{thm:triToFilling} by induction on dimension.  We first
define an action of $J^m(\R^k)$ on augmented maps.  Let $p=(p_1,p_2)\in
J^m(\R^k)$.  Define
$$p\cdot (\alpha, \{f_\Delta\}_{\Delta\subset X})=(p\cdot \alpha,
\{p\cdot f_\Delta\}_{\Delta\subset X})$$
where 
\begin{align*}
  (p\cdot \alpha)(x)&=p\alpha(x)\\
  (p\cdot f_\Delta)(x)&=f_\Delta(x-p_1)+P(p)(x).
\end{align*}
It is easily checked that this is a group action on the space of augmented maps.

\begin{lemma}\label{lemma:jetExist}
  There are $\tau, \eta, \phi, $ and $\psi$ satisfying the conditions of
  Theorem~\ref{thm:triToFilling} for $G=J^m(\R^{k+1})$ and $k=k$.
\end{lemma}
\begin{proof}
  Let $\Gamma=\Gamma_{m,k+1}$ and recall that $s_2(\Gamma)\subset \Gamma$.
  Let $p:J^m(\R^{k+1})\to \R^{k+1}$ be the bundle projection
  and let $(\tau,f:\tau\to G)$ be a $\Gamma$-equivariant
  triangulation of $G$ and let $(\eta,g:\eta\to G\times[0,1])$
  be an $s_2(\Gamma)$-equivariant triangulation of
  $G\times[0,1]$, as in Lemma~\ref{lem:trisExist}.  We will
  construct horizontal maps $\phi:\tau\to G$ and
  $\psi:\eta\to G$ by constructing augmented maps on the
  $k+1$-skeletons of $\tau$ and $\eta$. 

  We can construct an action of $\Gamma$ on $\R^{k+1}$ by letting the
  action of $\gamma$ send $x\mapsto x+p(\gamma)$.  After possibly
  subdividing $\tau$, we can construct a $\Gamma$-equivariant map
  $\kappa:\tau\to \R^{k+1}$ so that the vertices of any simplex of
  $\tau$ lie in general position and each simplex is mapped linearly
  to $\R^{k+1}$; this is an embedding on each simplex in
  $\tau^{(k+1)}$.
  
  We can now use Lemma \ref{lemma:simplexextension} to build a
  $\Gamma$-equivariant augmented map on the 0-skeleton, then the
  1-skeleton, and so on up to the $(k+1)$-skeleton.  This constructs a
  horizontal $\Gamma$-equivariant map on $\tau^{(k+1)}$ which can
  then be extended to all of $\tau$.  

  We construct $\psi$ by a similar process.  The main difference is
  the starting point; if $\tau_0$ and $\tau_1$ are as in
  Lemma~\ref{lem:trisExist} and if $\iota_i:\tau_i\to \eta$, $i=0,1$
  are the inclusions of the $\tau_i$ into $\eta$, we define $\psi_0:\tau_0\cup \tau_1\to G$
  by   
  $$\psi_0|_{\iota_0(\tau_0)}=\phi\circ \iota_0^{-1}$$
  and
  $$\psi_0|_{\iota_1(\tau_1)}=s_2\circ \phi\circ \iota_1^{-1}.$$
  and extend this to a $s_2(\Gamma)$-equivariant $(k+1)$-horizontal
  map on $\eta$.
\end{proof}

We thus have
\begin{thm} \label{thm:eucdehn}
  $J^m(\R^k)$ satisfies the filling inequalities
  $$\FV^n(V)\prec V^{\frac{n}{n-1}} \qquad \text{for $2\le n \le k$.}$$
\end{thm}

In the case of the Heisenberg group $H_{2k+1}=J^1(\R^k)$, we have
\begin{cor}
  $H_{2k+1}$ satisfies the filling inequalities
  $$\FV^n(V)\prec V^{\frac{n}{n-1}} \qquad \text{for $2\le n \le k$.}$$
\end{cor}

\section{More upper bounds}\label{sec:moreUpper}

When $G$ has no horizontal $(k+1)$-manifolds, the constructions of
Section~\ref{sec:fillApprox} still provide fillings of cycles, but
these fillings satisfy weaker bounds.  In some cases, however
(especially when there are many horizontal $k$-manifolds), these
weaker bounds may still be sharp.  One notable application is a sharp upper
bound for the filling volume of the higher-dimensional Heisenberg
groups in the middle dimension.

\begin{lemma}\label{lem:constRi2}
  Let $k>0$.  Let $\tau$, $\eta$,
  $\phi:\tau\to G$, and $\psi:\eta\to G$ be as in
  Lemma~\ref{lem:constRi}, except with no requirement that $\phi$ and
  $\psi$ be horizontal.

  Let $f:\R^+\to \R^+$ be a function such that for any $(k+1)$-simplex
  $\Delta$ of $\eta$, we have $\mass s_t(\psi(\Delta))\le f(t)$.  If,
  as in Lemma~\ref{lem:constRi}, we define
  $$R_i(\alpha):=s_{2^i}( R_{\psi(\eta)}(s_{2^{-i}}(\alpha))),$$ 
  there is a $c$ such that for all $i\ge 0$ and for all $k$-cycles $\alpha$, we have
  $$\partial R_i(\alpha)=P_{i+1}(\alpha)-P_{i}(\alpha)$$
  and
  $$\mass R_i(\alpha) \le  2^{-ki} c f(2^i) \mass \alpha.$$
\end{lemma}
\begin{proof}
  The fact that
  $$\partial R_i(\alpha)=P_{i+1}(\alpha)-P_{i}(\alpha)$$
  follows from the proof of Lemma~\ref{lem:constRi}.

  If $\beta$ is a simplicial chain in a simplicial complex $\sigma$,
  it can be written as $\beta=\sum_{\Delta\in\sigma} b_\Delta \Delta$,
  where $\Delta$ ranges over the simplices of $\sigma$.  Define 
  $$\|\beta\|_1=\sum_{\Delta\in\sigma} |b_\Delta|.$$

  Recall that if $\gamma$ is a Lipschitz $k$-cycle, then we defined
  $$R_{\psi(\eta)}(\gamma)=\psi_\sharp(P_\eta(X(\gamma))).$$
  By \ref{thm:fedflem}, there is a $c$ depending only on $\eta$ such that
  $\|P_\eta(X(\gamma))\|_1\le c\mass \gamma$.  Therefore, we have
  $$\|P_\eta(X(s_{2^{-i}}(\alpha)))\|_1\le c 2^{-ki}\mass \alpha,$$
  so $R_i(\alpha)$ is the sum of at most $(c 2^{-ki}\mass \alpha)$
  singular simplices of the form $(s_{2^i}\circ \psi)_\sharp(\Delta)$,
  where $\Delta$ is a simplex of $\eta$.  Each of these has mass $\le
  f(2^i)$, so 
  $$\mass R_i(\alpha) \le  2^{-ki} c f(2^i) \mass \alpha.$$
\end{proof}

The same construction used to prove Theorem~\ref{thm:triToFilling} can
be used to prove:
\begin{thm}\label{thm:triToFilling2}
  Let $k$, $\tau$, $\phi$, $\eta$, and $\psi$ satisfy the hypotheses
  of Lemma~\ref{lem:constRi2}.  Then
  $$\FV^{k+1}_G(V)\preceq \sum_{i=1}^{(\log_2 V)/k} f(2^i)
  \frac{V}{2^{ki}}.$$
  So if $f(t)\sim t^{p}$, where $p>k$, we have 
  $$\FV^{k+1}_G(V)\preceq V^{\frac{p}{k}}.$$
\end{thm}

To apply this to the jet groups, we use the maps constructed in
Lemma~\ref{lemma:jetExist}.  In that lemma, we constructed
triangulations $\tau$ and $\eta$ and equivariant $k$-horizontal maps
$\phi:\tau\to J^m(\R^{k})$ and $\psi:\eta\to J^m(\R^{k})$; these maps
are horizontal when restricted to $\tau^{(k)}$ and $\eta^{(k)}$, and
we will define equivariant maps $\phi'$ and $\psi'$ which extend
$\phi|_{\tau^{(k)}}$ and $\psi|_{\eta^{(k)}}$.

Let $\Delta_1,\dots,\Delta_d$ be a set of $(k+1)$-simplices which form
a fundamental domain for $\tau^{(k+1)}$.  Then $\phi'$ is already
defined on $\partial \Delta_i$ for each $i$ and is horizontal.  We can
identify $\Delta_i$ with the cone $C\partial\Delta_i:=(\partial
\Delta_i\times [0,1])/(\Delta_i\times \{0\})$ over $\partial \Delta_i$
and define $h_i:C\Delta_i\to G$ as the map which sends $(x,t)\in
\partial \Delta_i\times [0,1]$ to $s_t(\phi'(x))$.  This map is not
horizontal, but it scales relatively slowly; at almost every point in
the image, the tangent plane to $h_i$ is the sum of a $k$-dimensional
horizontal subspace and another vector.  Since $J^m(\R^k)$ has
nilpotency class $m+1$, any vector $v$ in $J^m(\R^k)$ satisfies
$\|s_t(v)\|\le t^{m+1}\|v\|$, and so we have
$$\mass (s_t\circ h_i)_\sharp(\Delta_i)\le t^{k+m+1}\mass
{h_i}_\sharp(\Delta_i).$$ Define $\phi'|_{\Delta_i}=h_i$ for each $i$,
and extend this definition to the entire $(k+1)$-skeleton by
equivariance.  Finally, extend $\phi'$ to an equivariant map on all of
$\tau$.  If $c=\max_i \vol h_{i}$, then this map has the property
that for any $(k+1)$-simplex $\Delta$ of $\tau$, we have $\mass
(s_t\circ \phi')_\sharp(\Delta)\le c t^{k+m+1}$.  

Similar techniques let us extend $\psi|_{\eta^{(k)}}$ to a map
$\psi':\eta\to G$ such that for any $(k+1)$-simplex $\Delta$ of
$\eta$, we have $\mass (s_t\circ \psi')_\sharp(\Delta)\le c'
t^{k+m+1}$.  We thus have

\begin{thm} \label{thm:middleJet}
  $J^m(\R^k)$ satisfies the filling inequality
  $$\FV^{k+1}(V)\prec V^{\frac{k+m+1}{k}}.$$
\end{thm}
\begin{cor} \label{thm:middleHeis}
  $H_{2k+1}$ satisfies the filling inequality
  $$\FV^{k+1}(V)\prec V^{\frac{k+2}{k}}.$$
\end{cor}

\section{Lower bounds}\label{sec:lower}

Gromov \cite{GroAII}, Baumslag, Miller, and Short\cite{BaMiSh}, Pittet
\cite{PittetCLB}, and Burillo \cite{Burillo} have all used
cohomological methods to prove lower bounds on filling functions in
nilpotent groups.  The idea is related to the use of calibrations in
geometric measure theory: if $\omega$ is a closed $k$-cochain in a group
$G$, then the value of $\omega$ on a $k$-chain $\alpha$ is determined by
its boundary $\partial{\alpha}$.  We denote this value by $\omega(\alpha)$.
If $\omega$ satisfies a bound $\omega(\beta)\le c \mass \beta$ for all
$k$-chains $\beta$ (for instance, if $\omega$ is a left-invariant de Rham
or simplicial cochain), it follows that $\FV^k(\partial \alpha)\ge
\omega(\alpha)/c.$

So one can find a lower bound on $\FV^k$ by finding a $k$-cocycle $\omega$
and a $k$-chain $\alpha$ such that $\partial\alpha$ is small and
$\omega(\alpha)$ is large.  Gromov gives one way of finding such chains in
the Heisenberg group, and we will sketch his argument below.  Let
$H_3=\R^3$ with the multiplication
$$(x,y,z)\cdot(x',y',z')=(x+x',y+y',z+z'+xy').$$
Then the 2-form $\omega=dy\wedge dz$ is closed and left-invariant, and
if $S$ is an open set in the $yz$-plane, then there is a $c\ne 0$ such
that $\int_S \omega=c\area S$.  In particular, if $S(y,z)$ is the
rectangle $\{0\}\times[0,y]\times [0,z]$, then $\int_{S(y,z)} z=cyz$,
so we can use these rectangles to produce 2-chains such that
$\omega(\alpha)$ is large.  It remains to make $\partial\alpha$ small.
Consider the closed curve $\gamma$ that connects $(0,0,0), (0,r,0),
(0,r,r^2),$ and $(0,0,r^2)$ by geodesics.  Since the $z$-axis of $H$
is quadratically distorted, this curve has length $\sim r$, but it is
close enough to the boundary of $S(r,r^2)$ that one can show that any
2-chain $\alpha$ which fills $\gamma$ has $\omega(\alpha)\sim r^3.$  This
shows that $\FV^2_{H_3}(r)\succeq r^3$.

Pittet \cite{PittetCLB} generalized this technique to other nilpotent groups by using a
2-dimensional subgroup in a Carnot group rather than the $yz$-plane in
the Heisenberg group.  He shows that if $G$ has a 2-dimensional
subgroup $A$ and if there is a closed invariant 2-form $\omega$ which
restricts to the area form on $A$, then one can construct cycles with
large filling areas by replacing the edges of large rectangles in $A$
by geodesics.  Burillo \cite{Burillo} generalized this argument to higher dimensions
and explicitly constructed families of cycles in the Heisenberg groups
to give lower bounds for their higher-order filling functions.

In the case that $G$ is a Carnot group, one can often use scaling
automorphisms to help construct difficult-to-fill cycles.  For
example, say that $\omega$ is a closed invariant $k$-form in $G$ with
the property that $s_t^*(\omega)=t^{d}\omega$.  In this case, if
$\alpha$ is a $k$-chain with horizontal boundary and if
$\omega(\alpha)\ne 0$, then we have $\mass\partial s_t(\alpha)\sim
t^{k-1}$, but 
$$\omega((s_t)_\sharp(\alpha))=[s_t^*(\omega)](\alpha)\sim t^d.$$
We will prove lower bounds for other nilpotent groups using this
construction.  We prove the following:
\begin{thm}\label{thm:lowerbounds}
  Let $G$ be a Carnot group and let $(\tau,f:\tau\to G)$ be a
  triangulation of $G$.  Let $\phi:\tau\to G$ be a
  $k$-horizontal map which is a bounded distance from $f$.
  Let $M$ be a $(k+1)$-dimensional subgroup of $G$.  If $\omega$ is a
  closed invariant $(k+1)$-form on $G$ which restricts to the volume
  form on $M$ and which satisfies $s_t^*(\omega)=t^{d}\omega$ for all
  $t>0$, then 
  $$\FV^{k+1}_G(V)\succ V^{d/k}$$
\end{thm}
\begin{proof}
  We want to find a $k+1$-chain $\alpha$ such that $\partial \alpha$
  is $k$-horizontal and $\omega(\alpha)\ne 0$.  We start with a large
  ball in $M$.  Let $\epsilon>0$ be a small number to be chosen later.
  Since any subgroup of a nilpotent group is nilpotent, $M$ has
  polynomial growth, so we can choose a ball $B$ such that $\vol \bd
  B<\epsilon \vol B$.  Let $[B]$ be the fundamental class of $B$,
  oriented so that $\mass [B]=\vol B=\omega([B])$.
  
  Consider the chain 
  $$\alpha=Q_{\phi(\tau)}(\partial [B])+[B].$$
  where $Q_{\phi(\tau)}$ is defined as in
  Section~\ref{sec:fillApprox}.  We have
  $$\partial\alpha=P_{\phi(\tau)}(\partial
  [B])-\partial[B]+\partial[B]=P_{\phi(\tau)}(\partial
  [B]),$$ which is horizontal.  It only remains to check that if
  $\epsilon$ is sufficiently small, then $\omega(\alpha)\ne 0$.  Let
  $c_Q$ be the constant from Lemma~\ref{lem:Qprops} and let
  $$\|\omega\|:=\sup \{\frac{|\omega(\beta)|}{\mass \beta}\},$$
  where $\beta$ ranges over all $(k+1)$-chains.  Since $\omega$ is
  left-invariant and thus bounded, this supremum exists.  If we take $\epsilon <
  c_Q^{-1}\| \omega\|^{-1}$, then
  \begin{align*}
    |\omega(Q_{\phi(\tau)}(\partial [B]))| & \le
    \|\omega\|c_Q \mass \partial [B] \\
    & \le \|\omega\|c_Q \epsilon \mass [B] < \mass [B]
  \end{align*}
  Thus 
  $$\omega(\alpha)=\omega([B])+\omega(Q_{\phi(\tau)}(\partial
  [B]))>0,$$
  as desired, and the bound on the filling volume follows by the
  argument above.
\end{proof}

We can use this to show that for the groups $J^m(\R^k)$, the upper
bounds on $\FV_{d}$ when $d\le k+1$ given by
Theorems~\ref{thm:eucdehn} and \ref{thm:middleJet} are sharp.  By
Lemma~\ref{lemma:jetExist}, there are a triangulation $(\tau,f:\tau\to
J^m(\R^k))$ and a $k$-horizontal map $\phi:\tau\to J^m(\R^k)$ that
satisfy the conditions of Theorem~\ref{thm:lowerbounds}.  Thus, when
$d\le k+1$, we can bound $\FV_{d}$ from below by finding an
$d$-dimensional subgroup $M_d\subset J^m(\R^k)$ and a closed invariant
$d$-form $\omega_d$ which restricts to the volume form of $M_d$.

The Lie algebra of $J^m(\R^k)$ can be decomposed as
$$\mathfrak{j}_{m,k}=\left(\R^k\oplus W_m\right)\oplus
W_{m-1}\oplus\dots\oplus W_0.$$
Let $x_i:\R^k\to \R$ be the coordinate projections, and let
$e_1,\dots,e_k\in \R^k$ be the standard basis.  Since $W_0$ is a
1-dimensional subspace, we can let $z:W_0\to \R$ be an isometry.  Let
$\mathfrak{m}_d=\langle e_1,\dots,e_d\rangle$ for $i\le k$ and
$\mathfrak{m}_{k+1}=\R^k\oplus W_0$.  These subalgebras correspond to
subgroups $M_1,\dots,M_{k+1}$ of $J^m(\R^k)$.

Let
$$\omega_d=dx_1\wedge \dots \wedge dx_d$$
for $d\le k$ and let
$$\omega_{k+1}=dx_1\wedge \dots \wedge dx_k \wedge dz.$$
Then $\omega_i$ is a closed, invariant $i$-form for $1\le i\le k+1$,
and $\omega_i$ restricts to the volume form on $M_i$.
Furthermore, $s_t^*(\omega_i)=t^{i}\omega_i$ when $i\le k$ and
$s_t^*(\omega_{k+1})=t^{k+m}\omega_i.$
We conclude that:
\begin{thm} \label{thm:jetLower}
  $J^m(\R^k)$ satisfies the filling inequalities:
  \begin{align*}
    \FV^{i}(V)\succeq V^{\frac{i}{i-1}}&\text{ for $i\le k$}\\
    \FV^{k+1}(V)\succeq V^{\frac{k+m}{k}}.&\\
  \end{align*}
\end{thm}
In the case that $m=1$ (so that $J^m(\R^k)$ is a higher-dimensional Heisenberg group), this
reduces to a result of Burillo \cite{BurLow}:
\begin{thm}[Burillo]
  $H_{2k+1}$ satisfies the filling inequalities:
  \begin{align*}
    \FV^{i}(V)\succeq V^{\frac{i}{i-1}}&\text{ for $i\le k$}\\
    \FV^{k+1}(V)\succeq V^{\frac{k+2}{k}}.&\\
  \end{align*}
\end{thm}

\def\cprime{$'$}
\providecommand{\bysame}{\leavevmode\hbox to3em{\hrulefill}\thinspace}
\providecommand{\href}[2]{#2}

\end{document}